\documentclass[journal]{IEEEtran}
\usepackage[utf8]{inputenc}
\usepackage[english]{babel}

\usepackage[colorlinks=true, allcolors=blue]{hyperref}

\usepackage{amsmath, amssymb}

\usepackage{algorithm, algpseudocode}

\newtheorem{theorem}{Theorem}
\newtheorem{lemma}{Lemma}
\newtheorem{assumption}{Assumption}
\newtheorem{remark}{Remark}
\newtheorem{definition}{Definition}
\newtheorem{proposition}{Proposition}

\newtheorem{example}{Example}

\usepackage{graphicx}
\graphicspath{ {figure/} }
\usepackage{subfigure}
\usepackage{tabularx}
\usepackage{threeparttable}
\usepackage{booktabs} 
\usepackage{multirow}

\usepackage[noadjust]{cite}
\usepackage{url}

\IEEEoverridecommandlockouts

\title{\LARGE \bf
Dynamical System Approach for Optimal Control Problems with Equilibrium Constraints Using Gap-Constraint-Based Reformulation}
\author{Kangyu Lin, \IEEEmembership{Graduate Student Member, IEEE}, and Toshiyuki Ohtsuka, \IEEEmembership{Senior Member, IEEE}
\thanks{This work was supported in part by the JSPS KAKENHI Grant Number JP22H01510 and JP23K22780. 
Kangyu Lin was supported by the CSC scholarship (No. 201906150138).
The authors are with the Systems Science Course, Graduate School of Informatics, Kyoto University, Kyoto, Japan.
Email: {\tt\small k-lin@sys.i.kyoto-u.ac.jp, ohtsuka@i.kyoto-u.ac.jp }.}
}

\begin{document}

\maketitle

\begin{abstract}
This study focuses on using direct methods (first-discretize-then-optimize) to solve optimal control problems for a class of nonsmooth dynamical systems governed by differential variational inequalities (DVI), called optimal control problems with equilibrium constraints (OCPEC). 
In the discretization step, we propose a class of novel approaches to smooth the DVI.
The generated smoothing approximations of DVI, referred to as gap-constraint-based reformulations, have computational advantages owing to their concise and semismoothly differentiable constraint system.
In the optimization step, we propose an efficient dynamical system approach to solve the discretized OCPEC, where a sequence of its smoothing approximations is solved approximately. 
This system approach involves a semismooth Newton flow, thereby achieving fast local exponential convergence.
We confirm the effectiveness of our method using a numerical example.


\end{abstract}

\section{Introduction}\label{section: introduction}

\subsection{Background, motivation and related works}\label{subsection: background, motivation and related works}

Recent advances have attempted to extend optimal control to the control tasks of nonsmooth dynamical systems (i.e., state or its time derivatives have discontinuities).
These tasks arise in several cutting-edge engineering problems ranging from robotics to autonomous driving  \cite{Kim2025contactimplicitmpc, nurkanovic2023phdthesis, Lin2024noninteriorpoint, stewart2010optimal, Vieira2019phdthesis, brogliato2020dynamical, lin2024gap}.
\textcolor{blue}{Differential variational inequalities (DVIs) \cite{pang2008differential}, a unified mathematical formalism for modeling nonsmooth systems}, have garnered significant attention owing to their ability to exploit system structures using the mature theory of variational inequalities (VIs) \cite{facchinei2003finite}.
\textcolor{blue}{This study considers optimal control problems (OCPs) for a class of nonsmooth systems governed by DVI, known as \textit{optimal control problems with equilibrium constraints} (OCPECs)}.

Direct methods (i.e., \textit{first-discretize-then-optimize}) are practical for solving OCP of smooth systems \cite{betts2010practical}. 
However, its extension to the OCPEC encounters great challenges:
In the discretization step, discretizing a DVI using time-stepping methods \cite{stewart2011dynamics} leads to incorrect sensitivities, which introduce spurious minima into the discretized OCPEC \cite{nurkanovic2023phdthesis};
In the optimization step, the discretized OCPEC is a difficult nonlinear programming (NLP) problem called \textit{mathematical programming with equilibrium constraints} (MPECs), which violates all constraint qualifications (CQs) required by NLP theories.
One approach to alleviating these difficulties is to smooth the DVI and then use the continuation method in the smoothing parameter.
However, the smoothed DVI behaves similarly to the nonsmooth system when the smoothing parameter is small, and the problems to be solved become increasingly difficult. 

\textcolor{blue}{This study aims to extend the applicability of direct methods to OCPEC. 
Thus, two critical problems need to be addressed:}
\begin{itemize}
    \item \textcolor{blue}{How can the DVI be smoothed to make the smoothing approximation of the discretized OCPEC easier to solve?}
    \item \textcolor{blue}{How can a sequence of smoothing approximations of the discretized OCPEC be solved efficiently?}
\end{itemize}

The smoothing of DVI is not straightforward because VI involves infinitely many inequalities.
Existing smoothing approaches replace the VI with its Karush--Kuhn--Tucker (KKT) conditions.
These approaches introduce Lagrange multipliers, thereby generating smoothing approximations with many additional constraints.
Our recent work \cite{lin2024successive} proposed a \textit{multiplier-free} smoothing approach, which generates a smaller smoothing approximation by using gap functions \cite{fukushima1992equivalent} to reformulate VI as a small number of inequalities.
A recent study \cite{yao2024overcoming} also used gap functions to reformulate bilevel programs.
However, these gap functions were shown to be only \textit{once continuously differentiable} when initially proposed.
Thus, solution methods presented in \cite{lin2024successive} and \cite{yao2024overcoming} only use the first-order derivatives of gap functions and achieve a slow local convergence rate. 

After smoothing the DVI, we can obtain the solution to the discretized OCPEC by solving a sequence of its smoothing approximations.
This is a methodology known as the \textit{continuation method} \cite{allgower2012numerical}, where the core idea is to solve a difficult problem by solving a sequence of easier subproblems.
Its standard implementation is to solve each subproblem \textit{exactly}.
However, the latter subproblems become increasingly difficult, thereby requiring more computational time.
An alternative implementation is to solve each subproblem \textit{approximately} while ensuring that the approximation error is bounded, or better yet, finally converges to zero.
This implementation can be regarded as a case of the \textit{dynamical system approach}, also known as the \textit{systems theory of algorithm} \cite{dorfler2024towards}, where the iterative algorithm is viewed as a dynamical system and studied from a system perspective.
Dynamical system approaches have a long history and remain vibrant in many real-world applications \cite{ohtsuka2004continuation, fazlyab2017prediction, allibhoy2025anytime, francca2023nonsmooth, belgioioso2024online}.

\subsection{Contribution and outline}\label{subsection: contribution and outline}

Our contributions are summarized as follows, which are our solutions to the problems listed in subsection \ref{subsection: background, motivation and related works}.

\begin{itemize}
    \item \textcolor{blue}{We propose a class of novel and general approaches using \textit{Auchmuty's gap functions} \cite{auchmuty1989variational} to smoothing the DVI.
    The proposed approach is \textit{multiplier-free} and thus generates a smaller smoothing approximation for DVI.
    Moreover, we strengthen the differentiability of gap functions from \textit{once continuous differentiability} to \textit{semismooth differentiability}, which allows us to exploit their second-order gradient information for locally fast-converging algorithms.}
    \item \textcolor{blue}{We propose a \textit{semismooth Newton flow} dynamical system approach to solve the discretized OCPEC and prove the local exponential convergence under standard assumptions (i.e., strict complementarity, constraint regularity, and positive definiteness of the reduced Hessian).
    The proposed dynamical system approach facilitates solving a difficult nonsmooth OCP efficiently by leveraging the mature theory and algorithm for smooth systems.}
\end{itemize}

The remainder of this paper is organized as follows. 
Section \ref{section: preliminaries and problem formulation} reviews background material and formulates the OCPEC;
Section \ref{section: proposed approaches to smoothing the DVI} presents a novel class of approaches to smoothing the DVI;
Section \ref{section: Dynamical system approach to solve OCPEC} presents an efficient dynamical system approach to solve a sequence of smoothing approximations of the discretized OCPEC;
Section \ref{section: numerical experiment} provides the numerical simulation;
and Section \ref{section: conclusion} concludes this study.

\section{Preliminaries and problem formulation}\label{section: preliminaries and problem formulation}

\subsection{Notation, nonsmooth analysis and variational inequalities}\label{subsection: notation, nonsmooth analysis and VI}

We denote the nonnegative orthant of $\mathbb{R}^n$ by $\mathbb{R}^n_{+}$.
Given a vector $ x \in \mathbb{R}^n$, we denote the Euclidean norm by $\| x \|_2 = \sqrt{x^T x}$, the open ball with center at $x$ and radius $r > 0$ by $\mathbb{B}(x, r) = \{y \in  \mathbb{R}^{n} \ | \ \|y - x \|_2 < r \}$, and the Euclidean projector of $x$ onto a closed convex set $K \subseteq \mathbb{R}^n$ by $\Pi_{K}(x) := \arg \min_{y \in K} \frac{1}{2} \|y - x \|_2^2$.
Given a differentiable function $f: \mathbb{R}^n \rightarrow \mathbb{R}^m$, we denote its Jacobian by $\nabla_x f \in \mathbb{R}^{m \times n}$.
We say that a function $f$ is $k$-th Lipschitz continuously differentiable ($LC^k$ in short) if its $k$-th derivative is Lipschitz continuous.

Let function $G: \Omega \rightarrow \mathbb{R}^m$ be locally Lipschitz continuous in an open set $\Omega \subseteq \mathbb{R}^n$.
Let $N_G$ be the set of points where $G$ is not differentiable.
The \textit{generalized Jacobian} of $G$ at $x \in \Omega$ is defined as $\partial G (x) = \text{conv} \ \{ H \in  \mathbb{R}^{m \times n} \ | \ H = \lim\limits_{k \rightarrow \infty}  \nabla_x G(x^k) \}$ with $\{ x^k \}_{k = 1}^{\infty} \rightarrow x$ and $ x^k \notin N_G$,
where $\text{conv} S$ is the convex hull of a set $S$.
We say that $\partial G (x)$ is \textit{nonsingular} if all matrices in $\partial G (x)$ are nonsingular.
We say that $G$ is \textit{semismooth} at $\Bar{x} \in \Omega$ if $G$ is also directionally differentiable\footnote{The directional derivative of $G$ at $\Bar{x}$ exists in all directions} at $\Bar{x}$ and $\lim\limits_{ x \rightarrow \bar{x}  } \frac{G(x) + H(\Bar{x} - x) - G(\Bar{x}) }{ \| x - \Bar{x}\|_2} = 0$ holds\footnote{This limit means that $\partial G $ provides a Newton approximation for $G$ at $\Bar{x}$.} for any $x$ in the neighborhood of $\Bar{x}$ and any $ H \in \partial G(x)$. 
We say that $\theta: \Omega \rightarrow \mathbb{R}$ with $\Omega \subseteq \mathbb{R}^n $ open, is \textit{semismoothly differentiable} ($SC^1$ in short) at $x \in \Omega$ if $\theta$ is $LC^1$ in a neighborhood of $x$ and $\nabla_x \theta$ is semismooth at $x$.
A vector-valued function is $SC^1$ if all its components are $SC^1$.      

Given a feasible set $\mathcal{C} := \{ x \in \mathbb{R}^{n_x} | \boldsymbol{h}(x) = 0, \ \boldsymbol{c}(x) \geq 0\}$, where $\boldsymbol{h}: \mathbb{R}^{n_x} \rightarrow \mathbb{R}^{n_h}$ and $\boldsymbol{c}: \mathbb{R}^{n_x} \rightarrow \mathbb{R}^{n_c}$ are continuously differentiable, 
let $\mathcal{I}(x^*) = \{i \in \{1, \cdots, n_c \} |  \boldsymbol{c}_i(x^*) = 0\}$ be the active set of a point $x^* \in \mathcal{C}$,
we say that \textit{linear independence CQ} (LICQ) holds at $x^* \in \mathcal{C}$ if vectors $\nabla_x \boldsymbol{h}_i(x^*)$ with $i \in \{1, \cdots, n_h \}$ and $\nabla_x \boldsymbol{c}_i(x^*)$ with $i \in \mathcal{I}(x^*)$ are linearly independent, and \textit{Mangasarian--Fromovitz CQ} (MFCQ)\footnote{Note that LICQ implies MFCQ.} holds at $x^* \in \mathcal{C}$ if vectors $\nabla_x \boldsymbol{h}_i(x^*)$ with $i \in \{1, \cdots, n_h \}$ are linearly independent and a vector $d_x \in \mathbb{R}^{n_x}$ exists such that $\nabla_x \boldsymbol{h}(x^*) d_x = 0$ and $\nabla_x \boldsymbol{c}_i(x^*) d_x > 0, \forall i \in \mathcal{I}(x^*)$.

Given a closed convex set $K \subseteq \mathbb{R}^{n_\lambda}$ and a continuous function $F: \mathbb{R}^{n_\lambda} \rightarrow \mathbb{R}^{n_\lambda}$, \textcolor{blue}{ the \textit{variational inequalities} \cite{facchinei2003finite}, denoted by VI$(K, F)$, is to find a vector $\lambda \in K$ such that $(\omega - \lambda)^T F(\lambda) \geq 0, \forall \omega \in K$}.
The solution set of VI$(K, F)$ is denoted by SOL$(K, F)$.
If $K$ is finitely representable, i.e., $K:= \{\lambda \in \mathbb{R}^{n_\lambda} | g(\lambda) \geq 0\}$ with $g: \mathbb{R}^{n_\lambda} \rightarrow \mathbb{R}^{n_g}$ a (concave) function,
then \textcolor{blue}{VI solutions can be represented in a finite form.
Specifically, if $\lambda \in$ SOL$(K, F)$ and MFCQ holds at $\lambda$, then there exist Lagrange multipliers $\zeta \in \mathbb{R}^{n_g}$ such that }
\begin{subequations}\label{equation: KKT-condition-based reformulation for VI}
    \begin{align}
        & F(\lambda) -  \nabla_{\lambda}g(\lambda)^T\zeta = 0, \\
        & 0 \leq \zeta \perp g(\lambda) \geq 0. \label{equation: KKT-condition-based reformulation for VI complementarity condition}
    \end{align}
\end{subequations}
\textcolor{blue}{  We refer to (\ref{equation: KKT-condition-based reformulation for VI}) as the \textit{KKT condition} of the VI$(K, F)$ }.

\subsection{Optimal control problem with equilibrium constraints}\label{subsection: optimal control problem with equilibrium constraints}
We consider the finite-horizon continuous-time OCPEC:
\begin{subequations}\label{equation: continuous time OCPEC}
    \begin{align}
        \min_{x(\cdot), u(\cdot),  \lambda(\cdot)} \ &  \int_0^T L_S(x(t), u(t), \lambda(t)) dt \\
        \text{s.t.} \  &  \Dot{x}(t) = f(x(t), u(t), \lambda(t)), \quad x(0) = x_0, \label{equation: continuous time OCPEC ODE} \\
                       &  \lambda(t) \in \text{SOL}(K, F(x(t), u(t), \lambda(t))), \label{equation: continuous time OCPEC VI} 
    \end{align}
\end{subequations}
with 
state $x: [0, T] \rightarrow \mathbb{R}^{n_x}$, 
control $u: [0, T] \rightarrow \mathbb{R}^{n_u}$, 
algebraic variable $\lambda: [0, T] \rightarrow \mathbb{R}^{n_\lambda}$, 
and stage cost $L_S: \mathbb{R}^{n_x} \times \mathbb{R}^{n_u} \times \mathbb{R}^{n_{\lambda}} \rightarrow \mathbb{R}$. 
We call (\ref{equation: continuous time OCPEC ODE}) (\ref{equation: continuous time OCPEC VI}) a DVI, with ordinary differential equation (ODE) function $f: \mathbb{R}^{n_x} \times \mathbb{R}^{n_u} \times \mathbb{R}^{n_\lambda} \rightarrow \mathbb{R}^{n_x}$, VI set $K \subseteq \mathbb{R}^{n_\lambda}$, and VI function $F: \mathbb{R}^{n_x} \times \mathbb{R}^{n_u} \times \mathbb{R}^{n_\lambda} \rightarrow \mathbb{R}^{n_\lambda}$. 
Note that $\lambda(t)$ does not exhibit any continuity properties and thereby introduces discontinuities in $x(t)$ and $\dot{x}(t)$.
We make the following assumption: 
\begin{assumption}\label{assumption: set and function of continuous time OCPEC}
    Set $K$ is closed, convex, finitely representable, and LICQ holds.
    Functions $L_S, f, F$ are $LC^2$. 
    \hfill \IEEEQEDopen
\end{assumption}

Solving continuous-time OCPs by direct (multiple shooting) methods \cite{betts2010practical} first requires discretizing the dynamical systems.
At present, the discretization of DVI is still based on the time-stepping method \cite{stewart2011dynamics}, which discretizes the ODE (\ref{equation: continuous time OCPEC ODE}) implicitly\footnote{\textcolor{blue}{For nonsmooth ODEs, the numerical integration method is required to be stiffly accurate, which prevents numerical chattering, and algebraically stable, which guarantees bounded numerical errors. One method that meets these requirements is the implicit Euler method. See subsection 8.4.1 in \cite{stewart2011dynamics}.} } and enforces the VI (\ref{equation: continuous time OCPEC VI}) at each time point $t_n \in [0, T]$.
This leads to an OCP-structured MPEC:
\begin{subequations}\label{equation: discretized time OCPEC}
    \begin{align}
    \min_{\boldsymbol{x}, \boldsymbol{u}, \boldsymbol{\lambda}} \ &  \sum^{N}_{n=1}L_{S}(x_n, u_n, \lambda_n) \Delta t, \\
     \text{s.t.}     \quad              & x_{n-1} + \mathcal{F}(x_n, u_n, \lambda_n) = 0, \label{equation: discretized time OCPEC ODE}\\
                                        & \lambda_n \in \text{SOL}(K, F(x_n, u_n, \lambda_n)),  \quad n = 1, \dots ,N,\label{equation: discretized time OCPEC VI}
    \end{align}
\end{subequations}
where
$x_n\in \mathbb{R}^{n_x}$ and $\lambda_n\in \mathbb{R}^{n_\lambda}$ are the values of $x(t)$ and $\lambda(t)$ at $t_n$,  
$u_n \in \mathbb{R}^{n_{u}}$ is the piecewise constant approximation of $u(t)$ in the interval $(t_{n-1}, t_n]$, 
$N$ is the number of stages, 
$\Delta t = T/N$ is the time step,
and $\mathcal{F}: \mathbb{R}^{n_x} \times \mathbb{R}^{n_u} \times \mathbb{R}^{n_\lambda} \rightarrow \mathbb{R}^{n_x}$ forms the implicit discretization of the ODE  (e.g., $\mathcal{F}(x_n, u_n, \lambda_n) = f(x_n, u_n, \lambda_n) \Delta t - x_n$) for implicit Euler method).
We define 
$\boldsymbol{x} = [x^T_1, \cdots, x^T_N]^T$,
$\boldsymbol{u}= [u^T_1, \cdots, u^T_N]^T$ 
and $\boldsymbol{\lambda}= [\lambda^T_1, \cdots, \lambda^T_N]^T$ 
to collect variables.

The numerical difficulties in solving (\ref{equation: discretized time OCPEC}) lie in two aspects:
First, \textcolor{blue}{in nonsmooth systems, the sensitivities of $x(t)$ w.r.t. parameters and variables (e.g., $x_0$ and controls) are \textit{discontinuous}} \cite{nurkanovic2023phdthesis}, which cannot be revealed by the numerical integration of $x(t)$ no matter how small $\Delta t$ we choose.
In other words, the gradient information of (\ref{equation: discretized time OCPEC}) does not match that of (\ref{equation: continuous time OCPEC}).
As a result, many spurious minima exist in (\ref{equation: discretized time OCPEC}).
\textcolor{blue}{Second, the equilibrium constraints (\ref{equation: discretized time OCPEC VI}) violate CQs at any feasible point}\footnote{\textcolor{blue}{In general, SOL$(K, F)$ is a discrete set, which leads to the CQ violation. However, CQs are satisfied in certain special cases, for example, SOL$(K, F)$ reduces to $F = 0$ when $K = \mathbb{R}^{n_\lambda}$. Such cases are not considered here.} }.
These difficulties prohibit us from using NLP solvers to solve (\ref{equation: discretized time OCPEC}), as the gradient-based optimizer will be trapped in spurious minima near the initial guess due to the wrong sensitivity or fail due to the lack of constraint regularity.

\textcolor{blue}{One approach to alleviating these difficulties is to smooth the DVI by relaxing the VI.
Some studies \cite{stewart2010optimal, nurkanovic2023phdthesis, Lin2024noninteriorpoint} revealed that the sensitivity of a smoothing approximation to the nonsmooth system is correct if the time step $\Delta t$ is sufficiently smaller than the smoothing parameter.
Moreover, the relaxation of VI can also recover the constraint regularity.
Existing smoothing approaches replace the VI (\ref{equation: discretized time OCPEC VI}) with its KKT condition (\ref{equation: KKT-condition-based reformulation for VI}) and further relax the complementarity condition (\ref{equation: KKT-condition-based reformulation for VI complementarity condition}) into a set of parameterized inequalities  \cite{hoheisel2013theoretical}.
These approaches introduce Lagrangian multipliers, thereby generating an NLP problem with numerous additional inequalities.
This motivates us to explore better smoothing approximations for DVI.}

\section{Proposed approaches to smoothing the DVI}\label{section: proposed approaches to smoothing the DVI}

\subsection{Gap-constraint-based reformulations}\label{subsection: gap-constraint-based reformulations}
Our smoothing approaches are based on two new VI reformulations. 
These reformulations are inspired by Auchmuty's study \cite{auchmuty1989variational} for solving VI$(K, F)$.
Since the function $F(x, u, \lambda)$ in (\ref{equation: discretized time OCPEC VI}) also includes variables $x,u$, we introduce an auxiliary variable $\eta = F(x, u, \lambda)$ to reduce the complexity\footnote{\textcolor{blue}{If the VI is simple (e.g., $F$ is affine), the use of $\eta$ can be avoided.}} and redefine Auchmuty's function and its variants \cite{auchmuty1989variational, facchinei2003finite} as follows.
\begin{definition}\label{definition: redefined Auchmuty's function and its variants}
    Let $K \subseteq \mathbb{R}^{n_{\lambda}}$ be a closed convex set 
    and $d: \mathbb{R}^{n_{\lambda}} \rightarrow \mathbb{R}$ be a strongly convex and $LC^3$ function.
    We define the following functions:
    \begin{itemize}
        \item \textit{Auchmuty's function} $L^c_{Au}: \mathbb{R}^{n_{\lambda}} \times \mathbb{R}^{n_{\lambda}} \times \mathbb{R}^{n_{\lambda}} \rightarrow \mathbb{R}$:
        \begin{equation*}
             L^c_{Au}(\lambda, \eta, \omega) = cd(\lambda) - cd(\omega) + (\eta^T - c\nabla_\lambda d(\lambda))(\lambda-\omega)
        \end{equation*}
        where $c > 0$ is a given constant.
        \item \textit{Generalized primal gap function} $\varphi^{c}_{Au}: \mathbb{R}^{n_{\lambda}} \times \mathbb{R}^{n_{\lambda}} \rightarrow \mathbb{R}$:
            \begin{equation}\label{equation: explicit expression of generalized primal gap constraint for OCPEC VI gap function}
                \varphi^{c}_{Au}(\lambda, \eta) = \sup_{\omega \in K} L^{c}_{Au}(\lambda, \eta, \omega) = L^{c}_{Au}(\lambda, \eta, \hat{\omega}^{c}),
            \end{equation}        
        where $\hat{\omega}^{c}$ is the solution to the parameterized problem:
        \begin{equation}\label{equation: explicit expression of the solution to generalized primal gap constraint}
           \hat{\omega}^c  = \omega^c(\lambda, \eta) = \arg \max_{\omega\in K} L^c_{Au}(\lambda, \eta, \omega).
        \end{equation}
        \item \textit{Generalized D-gap function} $\varphi^{ab}_{Au}: \mathbb{R}^{n_{\lambda}} \times \mathbb{R}^{n_{\lambda}} \rightarrow \mathbb{R}$:
        \begin{equation}\label{equation: explicit expression of generalized D gap constraint for OCPEC VI gap function}
            \varphi^{ab}_{Au}(\lambda, \eta) = \varphi^{a}_{Au}(\lambda, \eta) - \varphi^{b}_{Au}(\lambda, \eta),
        \end{equation}        
        with two constants $a$ and $b$ satisfying $b > a > 0$. 
        \hfill \IEEEQEDopen
    \end{itemize}
\end{definition}

The properties of $\varphi^{c}_{Au}$ and $\varphi^{ab}_{Au}$ are summarized below. 
\begin{proposition}\label{proposition: properties of gap function}
    \textcolor{blue}{The following three statements are valid for gap functions $\varphi^{c}_{Au}(\lambda, \eta)$ and $\varphi^{ab}_{Au}(\lambda, \eta)$:}
    \begin{itemize}
        \item \textcolor{blue}{$\varphi^{c}_{Au}(\lambda, \eta)$ and $\varphi^{ab}_{Au}(\lambda, \eta)$ are $SC^1$;}
        \item \textcolor{blue}{$\varphi^{c}_{Au}(\lambda, \eta) \geq 0, \forall \lambda \in K$, and $\varphi^{c}_{Au}(\lambda, \eta) = 0$ with $\lambda \in K $ if and only if $\lambda \in$ SOL$(K, \eta)$;}
        \item \textcolor{blue}{$\varphi^{ab}_{Au}(\lambda, \eta) \geq 0, \forall \lambda \in \mathbb{R}^{n_{\lambda}}$, and $\varphi^{ab}_{Au}(\lambda, \eta) = 0 $ if and only if $\lambda \in$ SOL$(K, \eta)$.} \hfill \IEEEQEDopen
    \end{itemize}
\end{proposition}
\begin{IEEEproof}
    \textcolor{blue}{Following from the differentiability of a function defined by the supremum (Th. 10.2.1, \cite{facchinei2003finite}), we can first write down the explicit formula for the gradient of $\varphi^c_{Au}$}, which is 
    $\nabla_{\lambda}\varphi^c_{Au}(\lambda, \eta) = \eta^T - c(\lambda - \hat{\omega}^c)^T \nabla_{\lambda \lambda}d(\lambda)$ and $\nabla_{\eta}\varphi^c_{Au}(\lambda, \eta) = (\lambda-c\hat{\omega}^c)^T$.
    \textcolor{blue}{Following from the differentiability of the solution to a parameterized convex minimization problem (Corollary 3.5, \cite{von2008sc}), we have that $\hat{\omega}^c = \omega^c(\lambda, \eta)$ is semismooth}, thereby $\varphi^c_{Au}$ is $SC^1$.
    The differentiability of $\varphi^{ab}_{Au}$ follows similarly because it is defined as the difference of two $\varphi^c_{Au}$.

    Regarding the second and third statements, their proofs can be carried out similarly to those of Theorems 10.2.3 and 10.3.3 in \cite{facchinei2003finite}.
    See the proof in Appendix \ref{subsection: proof of second and third statements in properties of gap function}.
\end{IEEEproof}

\begin{remark}
    \textcolor{blue}{Existing studies (Sections 10.2 and 10.3, \cite{auchmuty1989variational}) only establish the \textit{once continuous differentiability} of gap functions. 
    Here, we upgrade the differentiability to \textit{semismooth differentiability}. 
    This improvement plays an important role in the fast-convergent algorithms presented in Section \ref{section: Dynamical system approach to solve OCPEC}}. \hfill \IEEEQEDopen
\end{remark}

Inspired by these properties, we propose two new reformulations that transform the infinitely many inequalities defining the VI (\ref{equation: discretized time OCPEC VI}) into a small number of inequalities.
\begin{proposition}\label{proposition: generalized gap constraint based reformulation for OCPEC VI}
    $\lambda \in$ SOL$(K, F(x,u,\lambda))$ iff $(x, u, \lambda, \eta)$ satisfies a set of $n_{\lambda} $ equalities and $n_g + 1$ inequalities:
    \begin{subequations}\label{equation: generalized primal gap constraint based reformulation for OCPEC VI}
        \begin{align}
            & F(x, u, \lambda) - \eta = 0, \label{equation: generalized primal gap constraint based reformulation for OCPEC VI function} \\
            & g(\lambda) \geq 0, \label{equation: generalized primal gap constraint based reformulation for OCPEC VI set}\\
            & \varphi^{c}_{Au}(\lambda, \eta) \leq 0, \label{equation: generalized primal gap constraint based reformulation for OCPEC VI gap constraint}   
        \end{align}
    \end{subequations}    
    or a set of $n_{\lambda} $ equalities and one inequality:
    \begin{subequations}\label{equation: generalized D gap constraint based reformulation for OCPEC VI}
        \begin{align}
            & F(x, u, \lambda) - \eta = 0, \\
            & \varphi^{ab}_{Au}(\lambda, \eta) \leq 0. \label{equation: generalized D gap constraint based reformulation for OCPEC VI gap constraint}   
        \end{align}
    \end{subequations}  
    We refer to (\ref{equation: generalized primal gap constraint based reformulation for OCPEC VI}) and (\ref{equation: generalized D gap constraint based reformulation for OCPEC VI}) as the \textit {primal-gap-constraint-based reformulation} and \textit {D-gap-constraint-based reformulation} for VI$(K, F(x,u,\lambda))$, respectively. \hfill \IEEEQEDopen
\end{proposition}
\begin{IEEEproof}
    This is the direct result following from the second and third statement of Proposition \ref{proposition: properties of gap function}. 
\end{IEEEproof}

Proposition \ref{proposition: generalized gap constraint based reformulation for OCPEC VI} provides two new approaches to smooth the DVI.
Specifically, we replace the VI (\ref{equation: discretized time OCPEC VI}) with its reformulation (\ref{equation: generalized primal gap constraint based reformulation for OCPEC VI}) (resp. (\ref{equation: generalized D gap constraint based reformulation for OCPEC VI})) and then relax the gap constraint (\ref{equation: generalized primal gap constraint based reformulation for OCPEC VI gap constraint}) (resp. (\ref{equation: generalized D gap constraint based reformulation for OCPEC VI gap constraint})).
This leads to two parameterized NLP problems:
\begin{subequations}\label{equation: discretized time OCPEC relax generalized primal gap constraint}
    \begin{align}
    \mathcal{P}^{c}_{gap}(s): \quad & \min_{\boldsymbol{x}, \boldsymbol{u}, \boldsymbol{\lambda}} \sum^{N}_{n=1}L_{S}(x_n, u_n, \lambda_n) \Delta t,  \\
     \text{s.t.}     \quad              & x_{n-1} + \mathcal{F}(x_n, u_n, \lambda_n) = 0, \\
                                        & F(x_n, u_n, \lambda_n) - \eta_n = 0, \label{equation: discretized time OCPEC relax generalized primal gap constraint auxiliary variable}\\
                                        & g(\lambda_n) \geq 0, \label{equation: discretized time OCPEC relax generalized primal gap constraint VI set} \\ 
                                        & s - \varphi^{c}_{Au}(\lambda_n, \eta_n) \geq 0, \quad n = 1, \dots, N, \label{equation: discretized time OCPEC relax generalized primal gap constraint gap constraint}
    \end{align}
\end{subequations}
\vspace{-1.5em}
\begin{subequations}\label{equation: discretized time OCPEC relax generalized D gap constraint}
    \begin{align}
    \mathcal{P}^{ab}_{gap}(s): \quad & \min_{\boldsymbol{x}, \boldsymbol{u}, \boldsymbol{\lambda}} \sum^{N}_{n=1}L_{S}(x_n, u_n, \lambda_n) \Delta t, \\
     \text{s.t.}     \quad              & x_{n-1} + \mathcal{F}(x_n, u_n, \lambda_n) = 0, \\
                                        & F(x_n, u_n, \lambda_n) - \eta_n = 0, \label{equation: discretized time OCPEC relax generalized D gap constraint auxiliary variable}\\
                                        & s - \varphi^{ab}_{Au}(\lambda_n, \eta_n) \geq 0, \quad n = 1, \dots, N,\label{equation: discretized time OCPEC relax generalized D gap constraint gap constraint}
    \end{align}
\end{subequations}
where $s \geq 0$ is a scalar relaxation parameter.

We now summarize the favorable properties of the proposed reformulations (\ref{equation: discretized time OCPEC relax generalized primal gap constraint}) and (\ref{equation: discretized time OCPEC relax generalized D gap constraint}) for the discretized OCPEC (\ref{equation: discretized time OCPEC}).
\textcolor{blue}{First, they are \textit{multiplier-free} (i.e., establishing the equivalence\footnote{\textcolor{blue}{In other words, the proposed reformulations (\ref{equation: generalized primal gap constraint based reformulation for OCPEC VI}) and (\ref{equation: generalized D gap constraint based reformulation for OCPEC VI}) establish the equivalence with the VI from a primal perspective, while the KKT-condition-based reformulation (\ref{equation: KKT-condition-based reformulation for VI}) does so from a primal–dual perspective.}} without Lagrange multipliers and related constraints), thereby possessing a \textit{more concise constraint system}, as shown in the fourth column of Table \ref{tab: Comparison of different reformulations for discretized time OCPEC VI}}.
\textcolor{blue}{Second, they are \textit{semismoothly differentiable} regardless of the value of $s$ and thereby can be solved with any given $s$ using Newton-type methods.
The fifth column of Table \ref{tab: Comparison of different reformulations for discretized time OCPEC VI} compares the differentiability of various VI reformulations.}
Third, their feasible set is equivalent to that of the original problem (\ref{equation: discretized time OCPEC}) when $s = 0$ and \textit{exhibits a feasible interior} when $s > 0$ (Example \ref{example: complementarity constraints}).
Hence, although $\mathcal{P}^{c}_{gap}(s)$ and $\mathcal{P}^{ab}_{gap}(s)$ lack constraint regularity when $s = 0$ (Theorem \ref{theorem: CQ violation of gap constraint based reformulation}), their regularity is recovered when $s > 0$.
Thus, we can solve the original problem (\ref{equation: discretized time OCPEC}) using the continuation method that solves a sequence of $\mathcal{P}^{c}_{gap}(s)$ or $\mathcal{P}^{ab}_{gap}(s)$ with $s \rightarrow 0$.

\begin{table*}[!t] 
    \centering
    \caption{Comparison of different reformulation for the equilibrium constraints (\ref{equation: discretized time OCPEC VI})}
    \label{tab: Comparison of different reformulations for discretized time OCPEC VI}
    \begin{tabular}{ccccc}
        \toprule
            \textbf{Reformulation}     & \textbf{Relaxed constraints}    & \textbf{Relaxation strategy}              & \textbf{Sizes}   & \textbf{Differentiability} (under Assumption \ref{assumption: set and function of continuous time OCPEC}) \\
        \midrule
        \multirow{5}{*}{KKT-condition-based}  &  \multirow{5}{*}{complementarity constraint}  &  Scholtes (Sec. 3.1 \cite{hoheisel2013theoretical}),              &  $N(n_{\lambda} + 3n_g)$    & $LC^2$  \\
                                              &                                               &  Lin-Fukushima (Sec. 3.2 \cite{hoheisel2013theoretical})          &  $N(n_{\lambda} + 2n_g)$    & $LC^2$\\
                                              &                                               &  Kadrani (Sec. 3.3 \cite{hoheisel2013theoretical})                &  $N(n_{\lambda} + 3n_g)$    & $LC^2$\\
                                              &                                               &  Steffensen--Ulbrich (Sec. 3.4 \cite{hoheisel2013theoretical})    &  $N(n_{\lambda} + 3n_g)$    & twice continuously differentiable  \\ 
                                              &                                               &  Kanzow--Schwartz  (Sec. 3.5 \cite{hoheisel2013theoretical})      &  $N(n_{\lambda} + 3n_g)$    & once continuously differentiable  \\
        \midrule
        \multirow{2}{*}{Gap-constraint-based} & gap constraint (\ref{equation: generalized primal gap constraint based reformulation for OCPEC VI gap constraint}) & Generalized primal gap  &  $N(n_{\lambda} + n_g + 1)$ & $SC^1$ \\
                                              & gap constraint (\ref{equation: generalized D gap constraint based reformulation for OCPEC VI gap constraint}) &   Generalized D gap     &  $N(n_{\lambda} +1)$        & $SC^1$ \\
        \bottomrule
    \end{tabular}
\end{table*}

\begin{figure*}[!t]
  \centering
  \subfigure[Contour of $\varphi^c(\lambda, \eta)$.]{\includegraphics[width=1.75in]{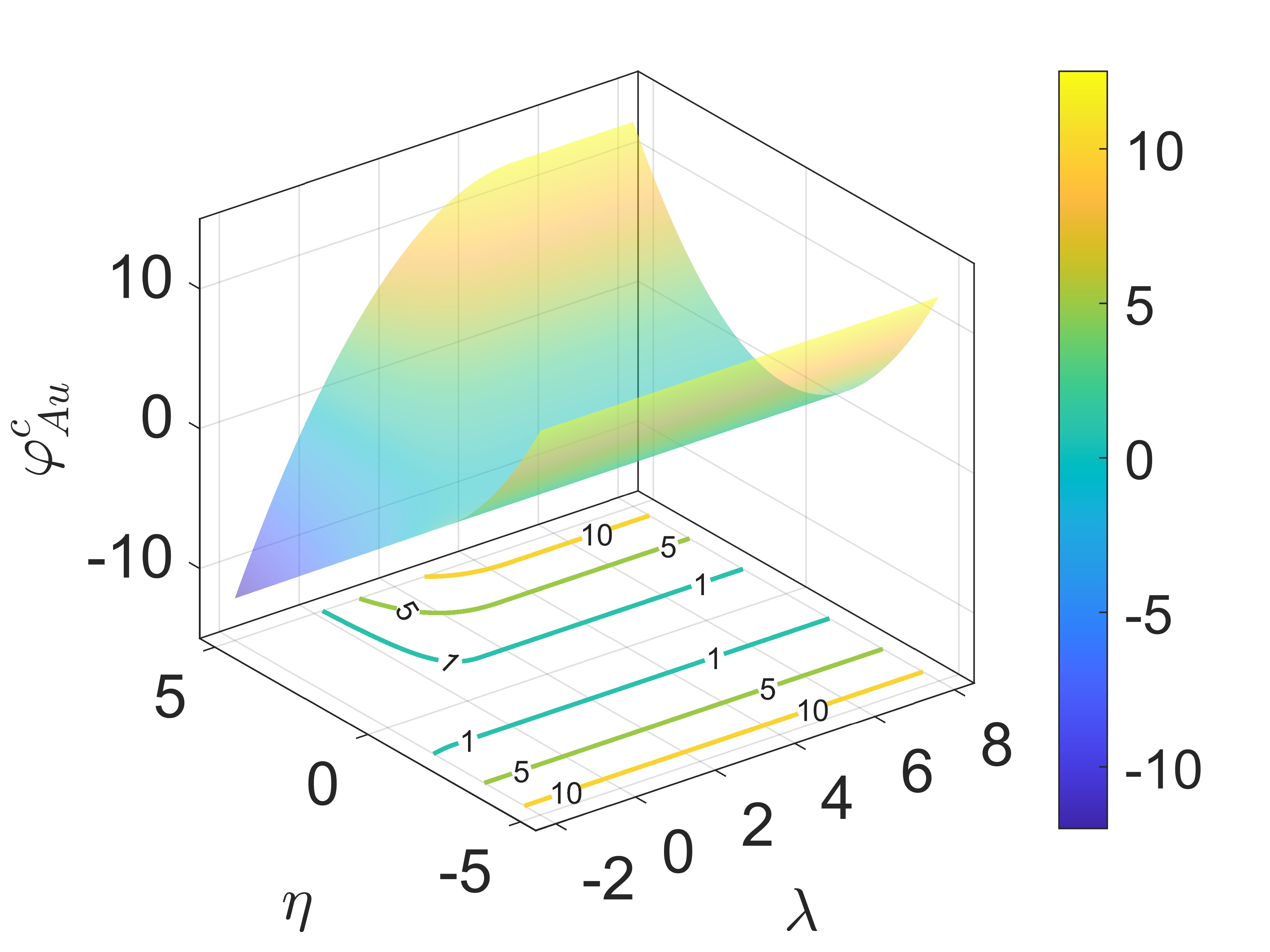}\label{fig: generalized primal gap function contour: complementarity constraint}}
  \hfil
   \subfigure[Feasible set of (\ref{equation: relax generalized primal gap constraint based reformulation for complementarity constraints}).]{\includegraphics[width=1.75in]{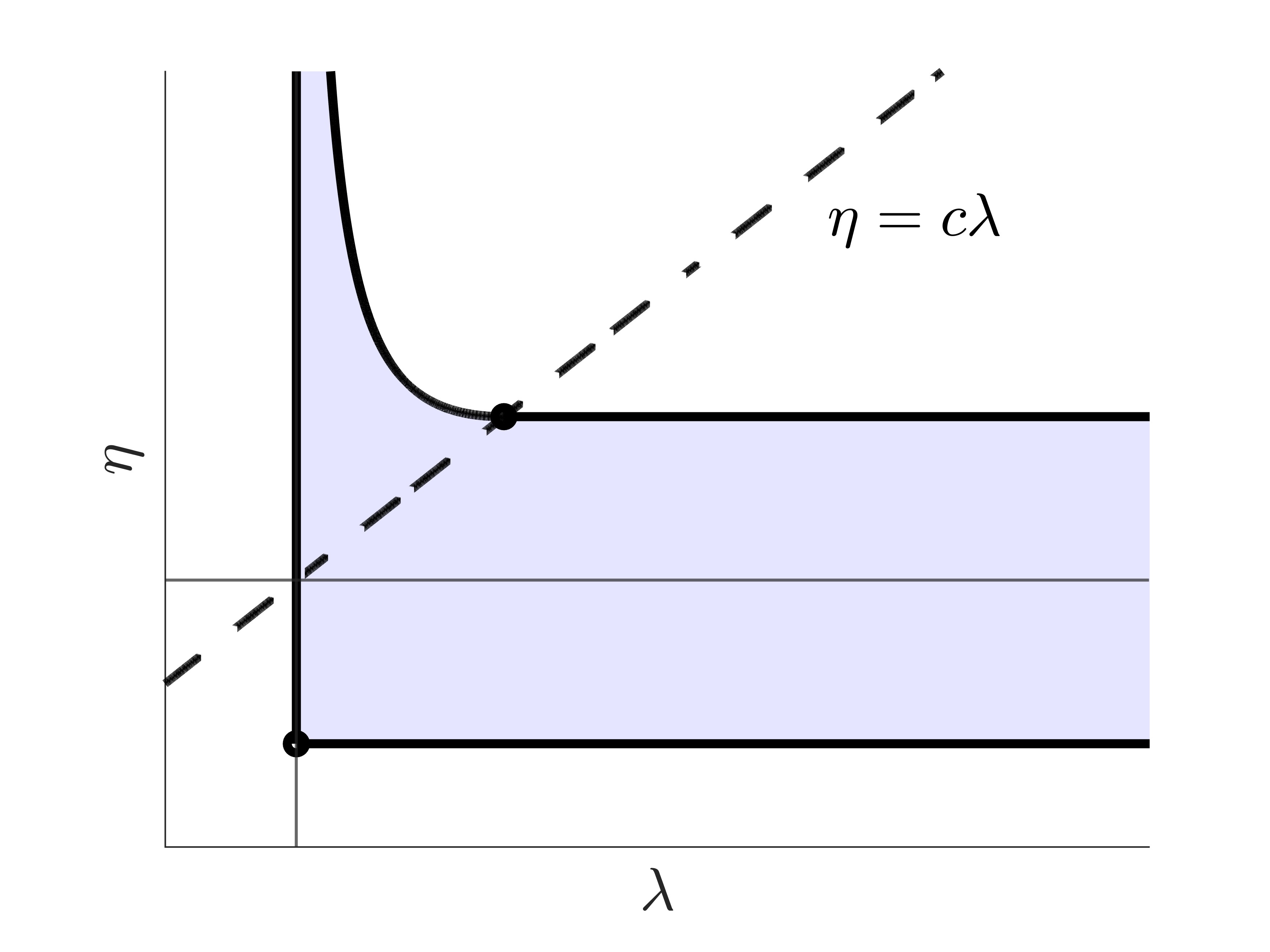}\label{fig: generalized primal gap function set: complementarity constraint}} 
  \hfil
  \subfigure[Contour of $\varphi^{ab}(\lambda, \eta)$.]{\includegraphics[width=1.75in]{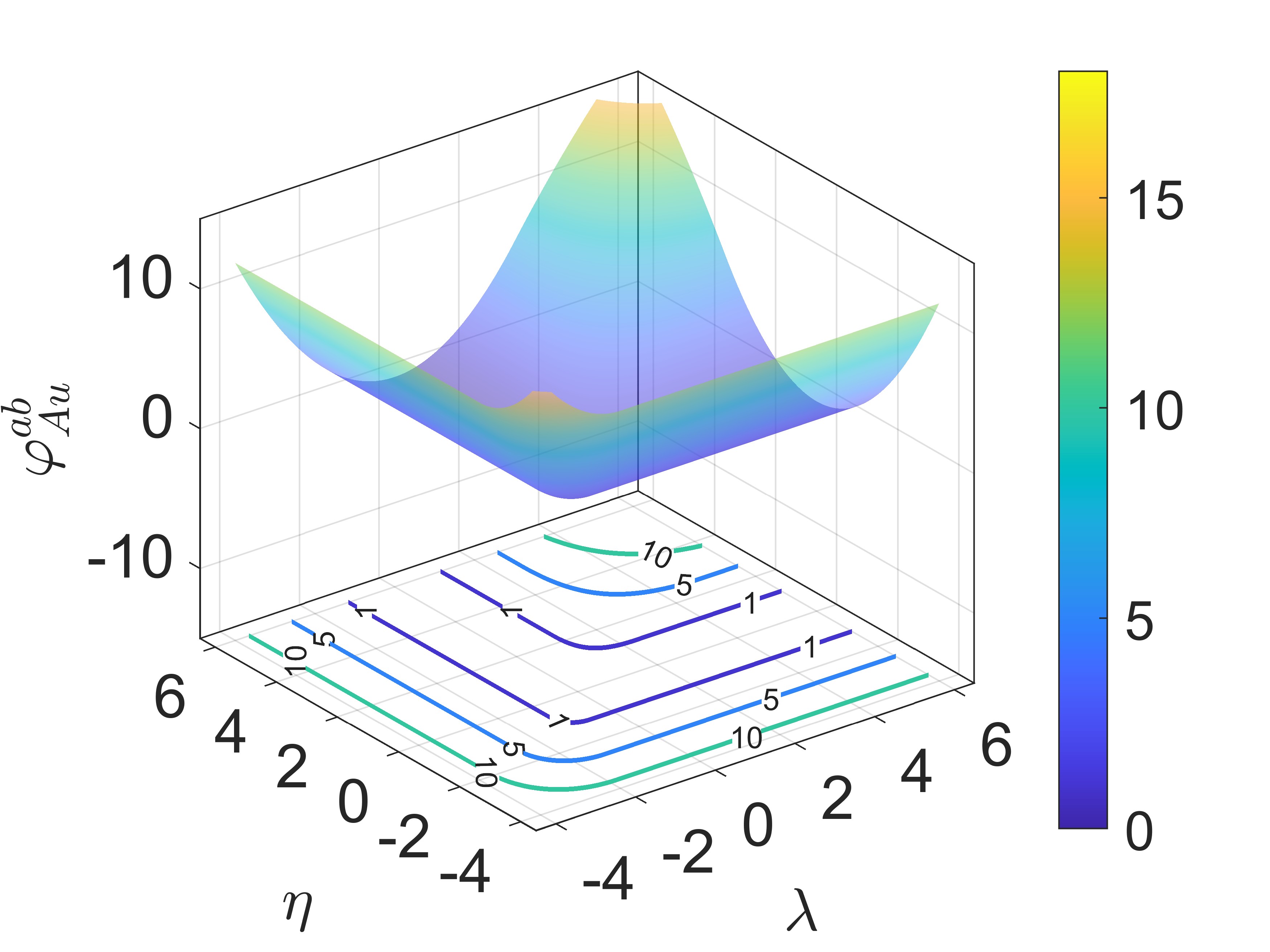}\label{fig: generalized D gap function contour: complementarity constraint}}
  \hfil
  \subfigure[Feasible set of (\ref{equation: relax generalized D gap constraint based reformulation for complementarity constraints}).]{\includegraphics[width=1.75in]{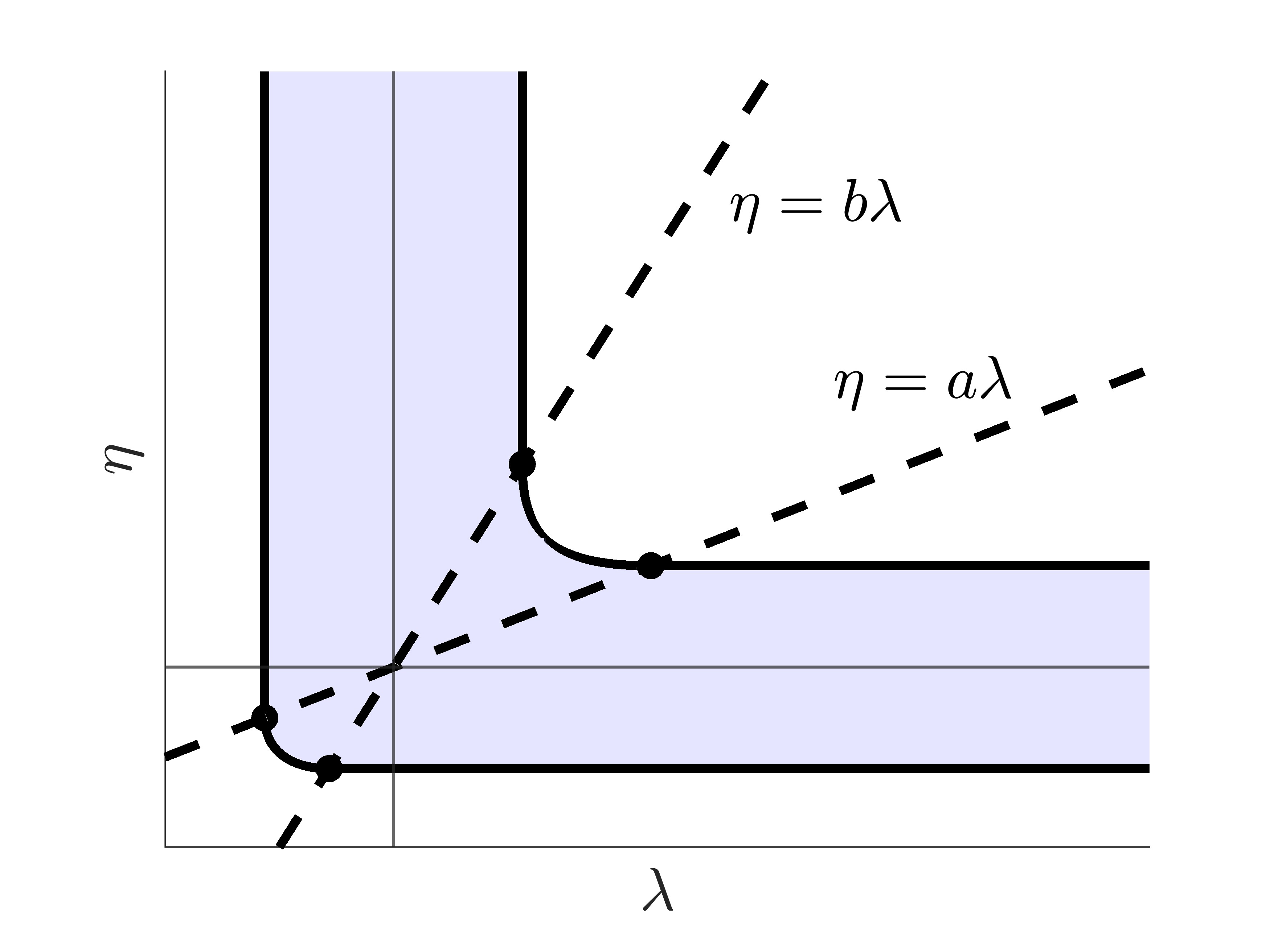}\label{fig: generalized D gap function set: complementarity constraint}}
  \caption{Geometric interpretation of the gap-constraint-based reformulations}
  \label{fig: geometric interpretation of the generalized gap constraint reformulation: complementarity constraint}
\end{figure*}

\subsection{Computation considerations, constraint regularity, and geometric interpretation}\label{subsection: computational considerations, constraint regularity and geometric interpretation}
Evaluating gap functions requires solving at least one constrained maximization problem, which typically is expensive.
Thus, we discuss how to exploit the OCP and VI structure to accelerate the evaluation of gap functions in (\ref{equation: discretized time OCPEC relax generalized primal gap constraint}) and (\ref{equation: discretized time OCPEC relax generalized D gap constraint}).

First, the maximization problems (\ref{equation: explicit expression of the solution to generalized primal gap constraint}) required to compute $\varphi^c_{Au}(\lambda_n, \eta_n)$ and $\varphi^{ab}_{Au}(\lambda_n, \eta_n)$ are \textit{stage-wise}, i.e., they involve only the variables and parameters of the same stage.
Moreover, the optimal active sets of the adjacent stage's problems may exhibit slight differences or even remain unchanged.
Thus, these problems can be solved in parallel with up to $N$ cores, or in serial using solvers with active set warm-start techniques \cite{lin2024successive}.
Second, the solution to problems (\ref{equation: explicit expression of the solution to generalized primal gap constraint}) may even possess an explicit expression.
For example, if $K$ is box-constrained (i.e., $K := \{\lambda \in \mathbb{R}^{n_\lambda} | b_l \leq \lambda \leq b_u \}$
with $b_l \in \{\mathbb{R} \cup \{-\infty\} \}^{n_{\lambda}}$ and $b_u \in \{\mathbb{R} \cup \{+\infty\} \}^{n_{\lambda}}$),
then we can specify $d(\cdot) = \frac{1}{2}\| \cdot \|^2_2$ to simplify (\ref{equation: explicit expression of the solution to generalized primal gap constraint}) as $\hat{\omega}^c = \omega^c(\lambda, \eta) = \Pi_{[b_l, b_u]}(\lambda - \frac{1}{c}\eta)$.
Moreover, the derivatives of $\Pi_{[b_l, b_u]}$, which are used to compute the second-order derivatives of gap functions, can be computed by the algorithmic differentiation software (e.g., CasADi \cite{Andersson2019}).

Next, we investigate whether the gap-constraint-based reformulations (\ref{equation: generalized primal gap constraint based reformulation for OCPEC VI}) and (\ref{equation: generalized D gap constraint based reformulation for OCPEC VI}) satisfy the constraint qualifications.
\begin{theorem}\label{theorem: CQ violation of gap constraint based reformulation}
    The gap-constraint-based reformulations (\ref{equation: generalized primal gap constraint based reformulation for OCPEC VI}) and (\ref{equation: generalized D gap constraint based reformulation for OCPEC VI}) violate the LICQ and MFCQ at any feasible point.  \hfill \IEEEQEDopen
\end{theorem}
\begin{IEEEproof}
    See the proof in Appendix \ref{subsection: proof of CQ violation of gap constraint based reformulation}.
\end{IEEEproof}

The violation of LICQ and MFCQ in the constraint systems (\ref{equation: generalized primal gap constraint based reformulation for OCPEC VI}) and (\ref{equation: generalized D gap constraint based reformulation for OCPEC VI}) is inevitable owing to their equivalences\footnote{\textcolor{blue}{Similar discussions also arise in bilevel optimization (Section IV-B in \cite{ouattara2018duality} and Section 4 in \cite{dempe2025duality}), where the CQs are interpreted as stating the constraints without the optima of an embedded optimization problem}, and the reformulations of the bilevel problem violate CQs once the equivalence holds.} with the VI solution set.
Nonetheless, the constraint systems (\ref{equation: generalized primal gap constraint based reformulation for OCPEC VI}) and (\ref{equation: generalized D gap constraint based reformulation for OCPEC VI}) have a feasible interior when their inequalities are relaxed.
Thus, if the constraint Jacobian of NLP problems (\ref{equation: discretized time OCPEC relax generalized primal gap constraint}) and (\ref{equation: discretized time OCPEC relax generalized D gap constraint}) satisfies certain full rank assumptions, then the LICQ and MFCQ hold on their constraint system when $s > 0$.

Finally, we provide a geometric interpretation of how the reformulations (\ref{equation: discretized time OCPEC relax generalized primal gap constraint}) and (\ref{equation: discretized time OCPEC relax generalized D gap constraint}) relax the equilibrium constraint (\ref{equation: discretized time OCPEC VI}) to smooth the DVI through a simple yet common example.
\begin{example}\label{example: complementarity constraints}
    Let $\lambda, \eta \in \mathbb{R}$ be scalar variables.
    We consider the complementarity constraint $0 \leq \lambda \perp \eta \geq 0$, which is a special case of VI$(K, \eta)$ with $K = \mathbb{R}_{+}$.
    The feasible set of $0 \leq \lambda \perp \eta \geq 0$ is the nonnegative part of axes $\lambda = 0$ and $\eta = 0$, which has no feasible interior point.
    By regarding $\lambda$ as the VI variable and specifying $d(\cdot) = \frac{1}{2}\| \cdot \|^2_2$, the reformulations (\ref{equation: generalized primal gap constraint based reformulation for OCPEC VI}) and (\ref{equation: generalized D gap constraint based reformulation for OCPEC VI}) for $0 \leq \lambda \perp \eta \geq 0$ are relaxed into:
    \begin{equation}\label{equation: relax generalized primal gap constraint based reformulation for complementarity constraints}
        \lambda \geq 0, \ s - \varphi^c_{Au}(\lambda, \eta) \geq 0,
    \end{equation}
    with $\varphi^c_{Au}(\lambda, \eta) =  \frac{1}{2c} \{\eta^2 - (\max(0, \eta - c \lambda))^2\}$, and
    \begin{equation}\label{equation: relax generalized D gap constraint based reformulation for complementarity constraints}
        s - \varphi^{ab}_{Au}(\lambda, \eta) \geq 0,
    \end{equation}
    with $\varphi^{ab}_{Au}(\lambda, \eta) = \varphi^{a}_{Au}(\lambda, \eta) - \varphi^{b}_{Au}(\lambda, \eta)$, respectively.
    The contour of $\varphi^c_{Au}$ and $\varphi^{ab}_{Au}$ are shown in Fig. \ref{fig: generalized primal gap function contour: complementarity constraint} and \ref{fig: generalized D gap function contour: complementarity constraint}.
    Hence, the feasible sets of (\ref{equation: relax generalized primal gap constraint based reformulation for complementarity constraints}) and (\ref{equation: relax generalized D gap constraint based reformulation for complementarity constraints}) are the colored regions in Fig. \ref{fig: generalized primal gap function set: complementarity constraint} and \ref{fig: generalized D gap function set: complementarity constraint}, which all exhibit a feasible interior.
    \hfill \IEEEQEDopen
\end{example}

\begin{remark}
    \textcolor{blue}{The choice of function $d$ mainly depends on the cost of computing $\hat{\omega}^c$ in (\ref{equation: explicit expression of the solution to generalized primal gap constraint}). 
    Under the requirements of Definition \ref{definition: redefined Auchmuty's function and its variants}, the simpler $d$ is, the better.}
    \textcolor{blue}{Parameters $a,b,c$ mainly affects the gradients of $\varphi^{c}_{Au}$ and $\varphi^{ab}_{Au}$.
    We recommend a moderate combination such as $a = 0.5, b = 2, c = 1$.}
    \hfill \IEEEQEDopen
\end{remark}

\section{Dynamical system approach to solve OCPEC}\label{section: Dynamical system approach to solve OCPEC}

\subsection{Problem setting and assumptions}\label{subsection: problem setting and assumptions}
\textcolor{blue}{At this stage, we can solve the discretized OCPEC (\ref{equation: discretized time OCPEC}) using the continuation method that solves a sequence of $\mathcal{P}^c_{gap}(s)$ or $\mathcal{P}^{ab}_{gap}(s)$ with $s \rightarrow 0$.
However, it is still difficult to solve $\mathcal{P}^c_{gap}(s)$ and $\mathcal{P}^{ab}_{gap}(s)$ when $s$ is small.
Thus, instead of solving each subproblem exactly using NLP solvers, we propose a novel dynamical system approach to perform the continuation method, which achieves a fast local convergence by exploiting the semismooth differentiability of the gap functions.}

Since both $\mathcal{P}^c_{gap}(s)$ and $\mathcal{P}^{ab}_{gap}(s)$ are parameterized NLPs, we consider the following NLP with parameterized inequalities throughout this section to stream the presentation:
\begin{subequations}\label{equation: general parameterized NLP}
    \begin{align}
        \mathcal{P}(s): \quad \min_{\boldsymbol{z}} \ & \boldsymbol{J}(\boldsymbol{z}), \\
        \text{s.t.}     \quad   & \boldsymbol{h}(\boldsymbol{z}) = 0, \label{equation: general parameterized NLP equality constraint}\\
                                & \boldsymbol{c}(\boldsymbol{z}, s) \geq 0, \label{equation: general parameterized NLP inequality constraint}
    \end{align}
\end{subequations}
with the decision variable $\boldsymbol{z} \in \mathbb{R}^{n_z}$, scalar parameter $s \geq 0$, and functions $\boldsymbol{J}: \mathbb{R}^{n_z} \rightarrow \mathbb{R}$, $\boldsymbol{h}: \mathbb{R}^{n_z} \rightarrow \mathbb{R}^{n_h}$ and $\boldsymbol{c}: \mathbb{R}^{n_z} \times \mathbb{R} \rightarrow \mathbb{R}^{n_c}$.
A point $\boldsymbol{z}$ satisfying (\ref{equation: general parameterized NLP equality constraint}) and (\ref{equation: general parameterized NLP inequality constraint}) is referred to as a \textit{feasible point} of $\mathcal{P}(s)$.
Let $\boldsymbol{\gamma}_h \in \mathbb{R}^{n_h}$ and $\boldsymbol{\gamma}_c \in \mathbb{R}^{n_c}$ be Lagrange multipliers, the KKT conditions for $\mathcal{P}(s)$ are:
\begin{subequations}\label{equation: KKT condition of the general parameterized NLP}
    \begin{align}
        & \nabla_{\boldsymbol{z}} \mathcal{L}(\boldsymbol{z}, \boldsymbol{\gamma}_h, \boldsymbol{\gamma}_c, s) = 0, \\
        & \boldsymbol{h}(\boldsymbol{z}) = 0, \\
        & 0 \leq \boldsymbol{c}(\boldsymbol{z}, s) \perp \boldsymbol{\gamma}_c \geq 0,  \label{equation: KKT condition of the general parameterized NLP (complementary condition)}      
    \end{align}
\end{subequations}
with Lagrangian $\mathcal{L}(\boldsymbol{z}, \boldsymbol{\gamma}_h, \boldsymbol{\gamma}_c, s) = \boldsymbol{J} + \boldsymbol{\gamma}_h^T\boldsymbol{h} - \boldsymbol{\gamma}^T_c \boldsymbol{c}$.
A triple $(\boldsymbol{z}^*, \boldsymbol{\gamma}_h^*, \boldsymbol{\gamma}_c^*)$ satisfying (\ref{equation: KKT condition of the general parameterized NLP}) is referred to as a \textit{KKT point} of $\mathcal{P}(s)$.
We make the following assumptions on $\mathcal{P}(s)$: 
\begin{assumption}\label{assumption: functions of general parameterized NLP}
    $\boldsymbol{J}$ and $\boldsymbol{h}$ are $LC^2$, whereas $\boldsymbol{c}$ is $SC^1$ w.r.t. $\boldsymbol{z}$ and affine in $s$;
\end{assumption}
\begin{assumption}\label{assumption: CQ violation of general parameterized NLP}
    Any feasible point violates MFCQ if $s = 0$;
\end{assumption}
\begin{assumption}\label{assumption: CQ recover of general parameterized NLP}
    If $s > 0$, then there exist at least one KKT point that satisfies the strict complementarity condition (i.e., $\boldsymbol{c}_i + \boldsymbol{\gamma}_{c,i} > 0, \forall i \in \{1, \cdots, n_c \}$) and LICQ; 
\end{assumption}
\begin{assumption}\label{assumption: reduced Hessian of general parameterized NLP}
    Let the generalized Jacobian of $\nabla_{\boldsymbol{z}} \mathcal{L}$ w.r.t. $\boldsymbol{z}$ be $\partial \nabla_{\boldsymbol{z}} \mathcal{L} \subset \mathbb{R}^{n_z \times n_z}$ and the elements of $\nabla_{\boldsymbol{z}} \mathcal{L}$ be $\nabla_{\boldsymbol{z}} \mathcal{L}_i$.   
    For each $\mathcal{H} \in \partial \nabla_{\boldsymbol{z}} \mathcal{L}$, we assume that $\mathcal{H} = \partial \nabla_{\boldsymbol{z}} \mathcal{L}_1 \times \partial \nabla_{\boldsymbol{z}} \mathcal{L}_2 \cdots \times  \partial\nabla_{\boldsymbol{z}} \mathcal{L}_{n_z}$\footnote{\textcolor{blue}{In general, $\mathcal{H} \in \partial \nabla_{\boldsymbol{z}} \mathcal{L}_1 \times  \cdots \times  \partial\nabla_{\boldsymbol{z}} \mathcal{L}_{n_z}$ (Proposition 7.1.14, \cite{facchinei2003finite}), and the inclusion is an equality when the non-differentiability of the components $\nabla_{\boldsymbol{z}}\mathcal{L}_{i}$ are unrelated (see Example 7.1.15, \cite{facchinei2003finite}).}} and the reduced Hessian $W^T \mathcal{H} W \succ 0$ at the KKT point, where $W \in \mathbb{R}^{n_z \times (n_z - n_h)}$ is a matrix whose columns are the basis for the null space of $\nabla_{\boldsymbol{z}}\boldsymbol{h}$.
    \hfill \IEEEQEDopen
\end{assumption}

Here, 
Assumption \ref{assumption: functions of general parameterized NLP} is consistent with Assumption \ref{assumption: set and function of continuous time OCPEC}, the differentiability of $\varphi_{Au}$ and $\varphi^{ab}_{Au}$, and the relaxation strategies (\ref{equation: discretized time OCPEC relax generalized primal gap constraint gap constraint}) and (\ref{equation: discretized time OCPEC relax generalized D gap constraint gap constraint}); 
Assumption \ref{assumption: CQ violation of general parameterized NLP} is consistent with Theorem \ref{theorem: CQ violation of gap constraint based reformulation},
and Assumptions \ref{assumption: CQ recover of general parameterized NLP} and \ref{assumption: reduced Hessian of general parameterized NLP} are used to ensure the nonsingularity of the KKT matrix, as shown in Lemma \ref{lemma: nonsingularity of the KKT matrix}.

\subsection{Fictitious-time semismooth Newton flow dynamical system}\label{subsection: Fictitious-time semismooth Newton flow dynamical system}

We now present the proposed dynamical approach to solve a sequence of $\mathcal{P}(s)$ with $s \rightarrow 0$.
We first transform the KKT system (\ref{equation: KKT condition of the general parameterized NLP}) into a system of semismooth equations.
This is achieved by using the Fisher-Burmeister (FB) function \cite{facchinei2003finite}: $\psi (a, b) = \sqrt{a^2 + b^2}-a-b$ with $a, b \in \mathbb{R}$. 
$\psi$ is semismooth and has the property that $\psi (a, b) = 0 \Leftrightarrow a \geq 0, b \geq 0, ab = 0$.
Let $\boldsymbol{v}_c \in \mathbb{R}^{n_c}$ be the auxiliary variable for $\boldsymbol{c}(\boldsymbol{z}, s)$.
We define $\boldsymbol{Y} = [\boldsymbol{z}^T, \boldsymbol{v}^T_c, \boldsymbol{\gamma}^T_h, \boldsymbol{\gamma}^T_c]^T \in \mathbb{R}^{n_Y}$ and rewrite (\ref{equation: KKT condition of the general parameterized NLP}) as\footnote{(\ref{equation: KKT condition of the general parameterized NLP (complementary condition)}) are mapped into $\Psi = 0$ using $\psi$ in an element-wise manner.}:
\begin{equation}\label{equation: KKT condition (equation system) of the general parameterized NLP}
    \boldsymbol{T}(\boldsymbol{Y}, s) = 
    \begin{bmatrix}
        \nabla_{\boldsymbol{z}} \mathcal{L}^T(\boldsymbol{z}, \boldsymbol{\gamma}_h, \boldsymbol{\gamma}_c, s)\\
        \boldsymbol{h}(\boldsymbol{z}) \\
        \boldsymbol{c}(\boldsymbol{z}, s) - \boldsymbol{v}_c \\
        \Psi(\boldsymbol{v}_c, \boldsymbol{\gamma}_c)
    \end{bmatrix}
    = 0,
\end{equation}
where the KKT function $\boldsymbol{T}: \mathbb{R}^{n_Y} \times \mathbb{R}_{+} \rightarrow \mathbb{R}^{n_Y}$ is semismooth based on Assumption \ref{assumption: functions of general parameterized NLP} and the semismoothness of $\psi$.

Let $\boldsymbol{Y}^*$ be a solution to (\ref{equation: KKT condition (equation system) of the general parameterized NLP}) with a given $s$.
We aim to find a solution $\boldsymbol{Y}^*$ associated with a small $s$.
Instead of considering $\boldsymbol{Y}^*$ as a function of $s$ and computing a sequence of solutions $\{\boldsymbol{Y}^{*,l} \}_{l=0}^{l_{max}}$ by solving (\ref{equation: KKT condition (equation system) of the general parameterized NLP}) exactly based on a given sequence of decreasing parameter $\{s^{l} \}_{l=0}^{l_{max}}$, we consider both $\boldsymbol{Y}^*$ and $s$ as functions of a \textit{fictitious time} $\tau \in [0, \infty)$, that is, we define the optimal solution trajectory and parameter trajectory as $\boldsymbol{Y}^*(\tau)$ and $s(\tau)$ respectively such that 
\begin{equation}\label{equation: continuous time KKT condition}
    \boldsymbol{T}(\boldsymbol{Y}^*(\tau), s(\tau))) = 0, \quad \forall \tau \geq 0.
\end{equation}

Regarding $s(\tau)$, since $s$ is a user-specified parameter, we define a dynamical system to govern $s(\tau)$:
\begin{equation}\label{equation: stable system to govern parameter trajectory}
    \dot{s} = - \epsilon_s (s - s_e), \ s(0) = s_0,
\end{equation}
where $\epsilon_s > 0$ is the stabilization parameter, and $s_0, s_e \in \mathbb{R}$ are the points where we expect $s(\tau)$ to start and converge. 

Regarding $\boldsymbol{Y}^*(\tau)$, let it start from $\boldsymbol{Y}^*(0) = \boldsymbol{Y}^*_0$, with $\boldsymbol{Y}^*_0$ a solution to (\ref{equation: KKT condition (equation system) of the general parameterized NLP}) associated with the given $s_0$.
Inspired by our earlier research in real-time optimization \cite{ohtsuka2004continuation}, we define a dynamical system evolving along the fictitious time axis such that its state $\boldsymbol{Y}(\tau)$, with $\boldsymbol{Y}(0) = \boldsymbol{Y}_0$ in the neighborhood of $\boldsymbol{Y}^*_0$, finally converge to $\boldsymbol{Y}^*(\tau)$ as $\tau \rightarrow \infty$.
\textcolor{blue}{This dynamical system is derived by stabilizing $\boldsymbol{T}(\boldsymbol{Y}(\tau), s(\tau)) = 0$ with a stabilization parameter $\epsilon_T > 0$:}
\begin{equation}\label{equation: differential equation for the KKT system (stabilized KKT residual form)}
    \dot{\boldsymbol{T}}(\boldsymbol{Y}(\tau), s(\tau)) = - \epsilon_T \boldsymbol{T}(\boldsymbol{Y}(\tau), s(\tau)),
\end{equation}
\textcolor{blue}{replacing the left-hand side of (\ref{equation: differential equation for the KKT system (stabilized KKT residual form)}) with the semismooth Newton approximation}\footnote{\textcolor{blue}{That is, $\dot{\boldsymbol{T}}(\boldsymbol{Y}(\tau), s(\tau)) = \mathcal{K} \dot{\boldsymbol{Y}} + \mathcal{S}  \dot{s}$.}} \textcolor{blue}{of $\boldsymbol{T}$, and substituting (\ref{equation: stable system to govern parameter trajectory}) into $\dot{s}$.
Consequently, we have:}
\begin{equation}\label{equation: stable system to govern solution trajectory}
    \dot{\boldsymbol{Y}} = -\mathcal{K}^{-1}(\epsilon_T \boldsymbol{T} - \epsilon_s \mathcal{S} (s - s_e)),
\end{equation}
\textcolor{blue}{with $\mathcal{K} \in \partial \boldsymbol{T} \subset \mathbb{R}^{n_Y \times n_Y}$ and $\mathcal{S} := \nabla_{s}\boldsymbol{T} \in \mathbb{R}^{n_Y}$.}
Here, $\partial \boldsymbol{T}$ is the generalized Jacobian of $\boldsymbol{T}$ w.r.t. $\boldsymbol{Y}$, and all KKT matrices $\mathcal{K}(\boldsymbol{Y}, s)$ have the form:
\begin{equation}\label{equation: KKT matrix}
    \mathcal{K} = 
    \begin{bmatrix}
        \mathcal{H} + \nu_H I   & 0 & \nabla_{\boldsymbol{z}} \boldsymbol{h}^T & -\nabla_{\boldsymbol{z}} \boldsymbol{c}^T \\
        \nabla_{\boldsymbol{z}} \boldsymbol{h} & 0 & -\nu_h I & 0 \\
        \nabla_{\boldsymbol{z}} \boldsymbol{c} & -I & 0 & 0 \\
        0  & \nabla_{\boldsymbol{v}_c} \Psi -\nu_c I & 0 & \nabla_{\boldsymbol{\gamma}_c} \Psi -\nu_c I 
    \end{bmatrix},
\end{equation}
where $\nu_H, \nu_h, \nu_c > 0$ are regularized parameters (e.g., $10^{-6}$) to ensure Assumptions \ref{assumption: CQ recover of general parameterized NLP} and \ref{assumption: reduced Hessian of general parameterized NLP}.
Matrix $\mathcal{S}$ is constant based on Assumption \ref{assumption: functions of general parameterized NLP}.
Finally, with the sampling of $s(\tau)$, we can compute $\boldsymbol{Y}(\tau)$ by numerically integrating (\ref{equation: stable system to govern solution trajectory})\footnote{\textcolor{blue}{Since $\epsilon_T$ and $\Delta \tau$ are user-specified, low-order integration schemes with lower computational complexity can still ensure accuracy and stability through appropriate choices of $\epsilon_T$ and $\Delta \tau$ (see Theorem \ref{theorem: error analysis of the Y system (explicit Euler)}).}}.  
In the following, we show that $\boldsymbol{Y}(\tau)$ converges to $\boldsymbol{Y}^*(\tau)$ exponentially.

\subsection{Convergence analysis}\label{subsection: convergence analysis}

First, we investigate the nonsingularity of KKT matrix.

\begin{lemma}\label{lemma: nonsingularity of the KKT matrix}
    For any given $s > 0$, let $\boldsymbol{Y}^*$ be the solution to (\ref{equation: KKT condition (equation system) of the general parameterized NLP}). 
    Every $\mathcal{K} \in \partial \boldsymbol{T}(\boldsymbol{Y}, s)$ is nonsingular for any $\boldsymbol{Y} \in \mathbb{R}^{n_Y}$ in the neighborhood of $\boldsymbol{Y}^*$.
    \hfill \IEEEQEDopen
\end{lemma}
\begin{IEEEproof}
    See the proof in Appendix \ref{subsection: proof of nonsingularity of the KKT matrix}.
\end{IEEEproof}

We now show the exponential convergence property.
\begin{theorem}\label{theorem: exponential convergence of the Y system}
    Let $\boldsymbol{Y}(\tau)$ and $s(\tau)$ be the trajectories governed by (\ref{equation: stable system to govern solution trajectory}) and (\ref{equation: stable system to govern parameter trajectory}), respectively.
    Let $\boldsymbol{Y}^*(\tau)$ be an optimal solution trajectory satisfying (\ref{equation: continuous time KKT condition}) and starting from $\boldsymbol{Y}^*(0) = \boldsymbol{Y}^*_0$, where $\boldsymbol{Y}^*_0$ is a solution to (\ref{equation: KKT condition (equation system) of the general parameterized NLP}) associated with the given $s_0$.
    Then, there exists a neighborhood of $\boldsymbol{Y}^*_0$ denoted by 
    $\mathcal{N}^{*}_{exp}$, such that for any $\boldsymbol{Y}(0) = \boldsymbol{Y}_0 \in \mathcal{N}^{*}_{exp}$, we have that $\boldsymbol{Y}(\tau)$ exponentially converges to $\boldsymbol{Y}^*(\tau)$ as $\tau \rightarrow \infty$, that is:
    \begin{equation}\label{equation: exponential convergence of the Y system}
        \|\boldsymbol{Y}(\tau) - \boldsymbol{Y}^*(\tau) \|_2 \leq k_1 \|\boldsymbol{Y}(0) - \boldsymbol{Y}^*(0) \|_2 e^{- k_2 \tau},
    \end{equation}
    with constants $k_1, k_2 > 0$.
    \hfill \IEEEQEDopen
\end{theorem}
\begin{IEEEproof}
    See the proof in Appendix \ref{subsection: proof of convergence and error analysis (explicit Euler) of the Y system}.
\end{IEEEproof}

\begin{remark}
    The exponential convergence of (\ref{equation: stable system to govern solution trajectory}) is a standard result if $\boldsymbol{T}$ is \textit{continuously differentiable}, which requires NLP functions in (\ref{equation: general parameterized NLP}) to be $LC^2$ (Proposition 2, \cite{fazlyab2017prediction}).
    Here we weaken the differentiability assumption by showing that the exponential convergence holds even if $\boldsymbol{T}$ is \textit{semismooth}, which only requires functions in (\ref{equation: general parameterized NLP}) to be $SC^1$.
    \hfill \IEEEQEDopen
\end{remark}

Finally, we provide an error analysis for the implementation of (\ref{equation: stable system to govern solution trajectory}) using the explicit Euler method.
\begin{theorem}\label{theorem: error analysis of the Y system (explicit Euler)}
    \textcolor{blue}{Let $\boldsymbol{Y}^*_l$, $\boldsymbol{Y}_l$ and $s_l$ be the points of trajectory $\boldsymbol{Y}^*(\tau)$, $\boldsymbol{Y}(\tau)$, $s(\tau)$ at $\tau = \tau_l$, respectively, and $ \dot{\boldsymbol{Y}}_l$ be the value of (\ref{equation: stable system to govern solution trajectory}) with $\boldsymbol{Y}_l$ and $s_l$.
    If $\{ \boldsymbol{Y}_l \}_{l=0}^{l_{max}}$ is updated by integrating (\ref{equation: stable system to govern solution trajectory}) using the explicit Euler method, i.e., $\boldsymbol{Y}_{l+1} = \boldsymbol{Y}_l + \Delta \tau \dot{\boldsymbol{Y}}_l$, then the following one-step error bound holds:}
    \begin{equation}
        \textcolor{blue}{\|\boldsymbol{Y}_{l+1} - \boldsymbol{Y}^*_{l+1} \|_2 \leq | 1 - \epsilon_T\Delta\tau| \|\boldsymbol{Y}_{l} - \boldsymbol{Y}^*_{l} \|_2 + \xi(s_l),}
    \end{equation}
    \textcolor{blue}{with $\xi(s_l) = k_3 \epsilon_s (s_l - s_e)$ and $k_3 > 0$ is a constant.}
    \hfill \IEEEQEDopen
\end{theorem}
\begin{IEEEproof}
    See the proof in Appendix \ref{subsection: proof of convergence and error analysis (explicit Euler) of the Y system}.
\end{IEEEproof}

\begin{remark}
    \textcolor{blue}{Theorem \ref{theorem: error analysis of the Y system (explicit Euler)} provides a criterion for error stability: $| 1 - \epsilon_T\Delta\tau| < 1$, which is consistent with the stability condition of the explicit Euler method.
    It also indicates that choosing a smaller $\epsilon_s$ can yield a tighter bound on the error. }
    \hfill \IEEEQEDopen
\end{remark}

\section{Numerical experiment}\label{section: numerical experiment}

The proposed method\footnote{\textcolor{blue}{The code is available at \url{https://github.com/KY-Lin22/Gap-OCPEC}.}} is implemented in MATLAB 2023b based on the CasADi symbolic framework \cite{Andersson2019}.
All experiments were performed on a laptop with a 1.80 GHz Intel Core i7-8550U.
We discretize the OCPEC (\ref{equation: continuous time OCPEC}) with $\Delta t = 5 \times10^{-4}$ into a parameterized NLP (\ref{equation: general parameterized NLP}) using gap-constraint-based reformulations (\ref{equation: discretized time OCPEC relax generalized primal gap constraint}) and (\ref{equation: discretized time OCPEC relax generalized D gap constraint}), where $\varphi^c_{Au}$ and $\varphi^{ab}_{Au}$ are specified with various $c$ and $a,b$, respectively.
\textcolor{blue}{We specify $\epsilon_s = 10$, $s_0 = 1$, $s_e = 10^{-3}$, $\epsilon_T = 50$ and compute $\boldsymbol{Y}(\tau)$ at each $\tau_l = l \Delta \tau$ by integrating $\dot{\boldsymbol{Y}}$ using the explicit Euler method.
We set the continuation step $l \in \{0, \cdots, l_{max} \}$ with $l_{max} = 500$, $\Delta \tau = 10^{-2}$, and obtain $\boldsymbol{Y}(0)$ by solving (\ref{equation: general parameterized NLP}) exactly with $s_0$ using a well-developed interior-point-method (IPM) NLP solver called IPOPT \cite{wachter2006implementation}}.

\begin{figure}
    \centering
    \includegraphics[width=0.9\linewidth]{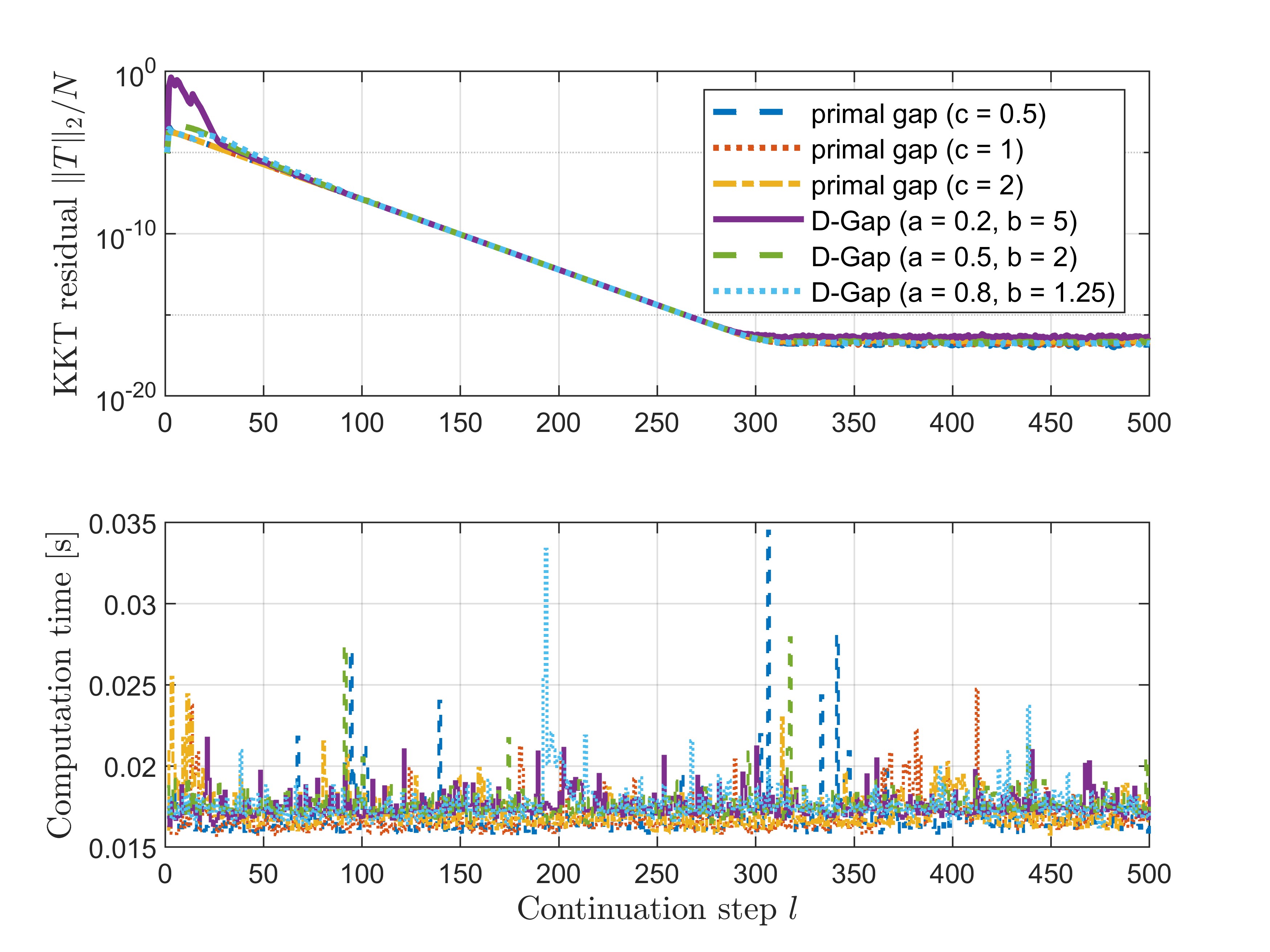}
    \caption{History of KKT residual and computation time (Proposed method).}
    \label{fig: KKT residual and computation time w.r.t the continuation step (Proposed method)}
\end{figure}

\begin{figure}
    \centering
    \includegraphics[width=0.9\linewidth]{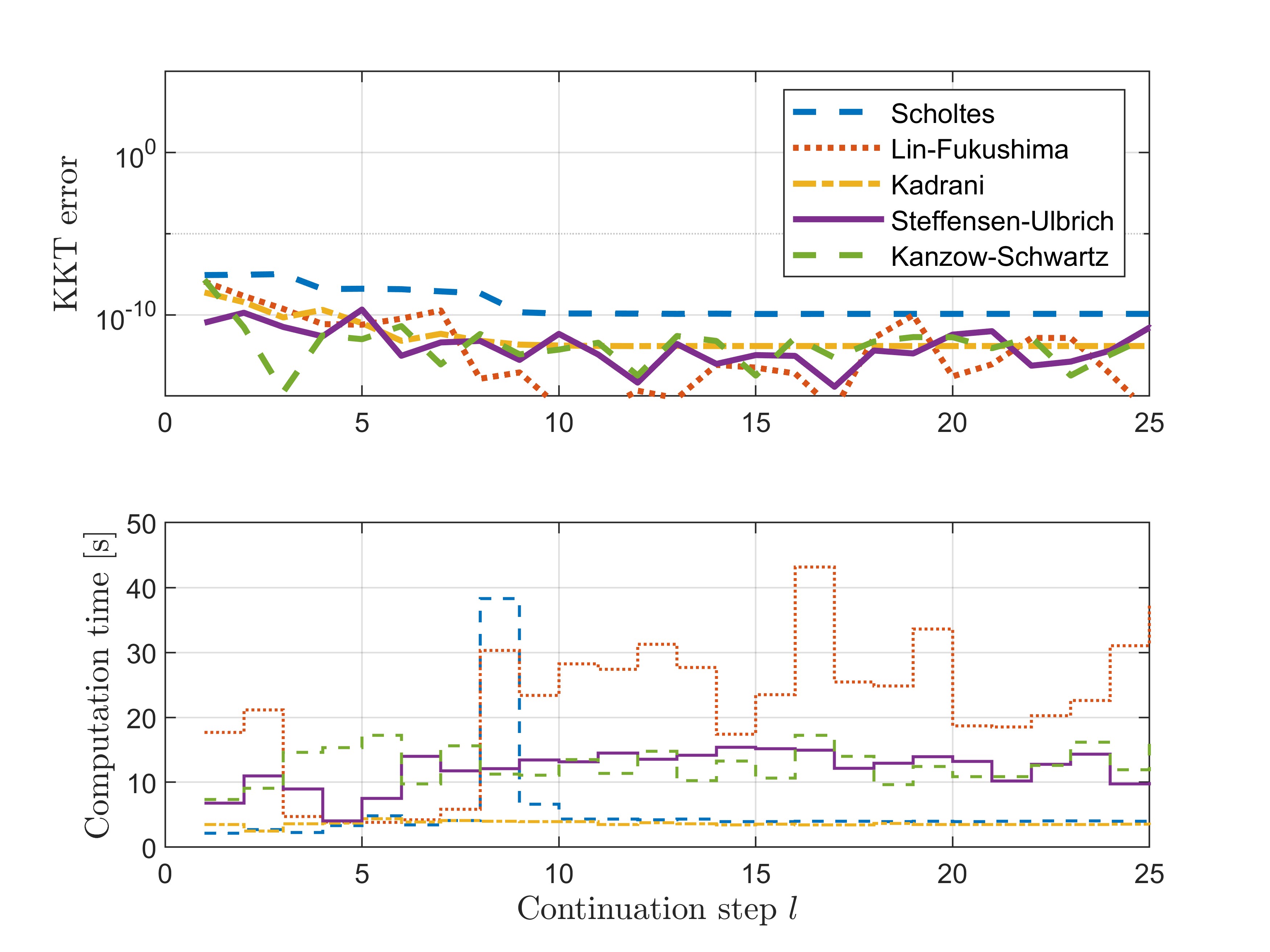}
    \caption{History of KKT error and computation time (Classical method).}
    \label{fig: KKT error and computation time w.r.t the continuation step (Classical methods)}
\end{figure}

The numerical example is an OCP of the linear complementarity system (LCS) taken from Example 7.1.5 of \cite{Vieira2019phdthesis}:
\begin{subequations}\label{equation: continuous OCP in LCS}
    \begin{align}
        \min_{x(\cdot), u(\cdot),  \lambda(\cdot)} \ &  \int_0^1 ( \| x(t) \|^2_2 + u(t)^2 +  \lambda(t)^2) dt \\        
        \text{s.t.} \  \dot{x}(t) =  &  \begin{bmatrix} 5 & -6 \\ 3 & 9 \end{bmatrix} x(t) + \begin{bmatrix} 0 \\ -4 \end{bmatrix} u(t) + \begin{bmatrix}4 \\ 5 \end{bmatrix} \lambda(t), \label{equation: continuous OCP in LCS (ODE)} \\
                       \eta(t) = &  \begin{bmatrix} -1 \\ 5 \end{bmatrix} x(t) + 6 u(t) +  \lambda(t),  \\
                        0 \leq & \ \lambda(t) \perp \eta(t) \geq 0, \label{equation: continuous OCP in LCS (complementarity constraints)} 
    \end{align}
\end{subequations}
with $x(0) = [-0.5, -1]^T$.
LCS is a special case of the DVI with affine functions $F, f$ and the VI set $K = \mathbb{R}_{+}$.
Thus, the gap functions $\varphi^c_{Au}$ and $\varphi^{ab}_{Au}$ for (\ref{equation: continuous OCP in LCS (complementarity constraints)}) have explicit expressions as discussed in Subsection \ref{subsection: computational considerations, constraint regularity and geometric interpretation}.
We solve this OCP using the proposed reformulations and dynamical system approach.
The history of the scaled KKT residual $\|\boldsymbol{T} \|_2/N$ w.r.t. the continuation step is shown at the top of Fig. \ref{fig: KKT residual and computation time w.r.t the continuation step (Proposed method)}.
The plots are linear on the log scale before converging to the point within machine accuracy.
Thus, the local exponential convergence is confirmed.
Moreover, as shown at the bottom of Fig. \ref{fig: KKT residual and computation time w.r.t the continuation step (Proposed method)}, the computation time for each continuation step remains nearly constant, ranging from 0.015 s to 0.035 s, with most values below 0.020 s.
\textcolor{blue}{For comparison, we also solve this example using classical methods, in which (\ref{equation: continuous OCP in LCS (complementarity constraints)}) is relaxed using the strategies presented in \cite{hoheisel2013theoretical}, and each subproblem is then solved exactly by IPOPT (the standard implementation of the continuation method).
The classical methods use a relaxation parameter sequence generated by discretizing (\ref{equation: stable system to govern parameter trajectory}) with the RK4 method using $\Delta \tau = 0.2$ ($\epsilon_s, s_0, s_e$ are the same as those used in the proposed method)}.
As shown in Fig. \ref{fig: KKT error and computation time w.r.t the continuation step (Classical methods)}, although the IPM KKT error of each continuation step remains within the desired small tolerance, each step requires a large amount of computation time, ranging from 1 s to 45 s.
This demonstrates the computational advantages of the proposed method.


\section{Conclusion}\label{section: conclusion}
This study focused on using the direct method to solve the OCPEC.
We addressed the numerical difficulties by proposing a new approach to smoothing the DVI and a dynamical system approach to solve a sequence of smoothing approximations of the discretized OCPEC.
The fast local convergence properties and computational efficiency were confirmed using a numerical example.
Our future work mainly focuses on incorporating a more sophisticated feedback structure into the dynamical system approach to achieve global convergence.


\appendix

\subsection{Proof of Theorem \ref{theorem: CQ violation of gap constraint based reformulation}}\label{subsection: proof of CQ violation of gap constraint based reformulation}

\begin{IEEEproof}
    Regarding the LICQ, Proposition \ref{proposition: properties of gap function} implies that the zeros of $\varphi^c_{Au}$ within the set $K$ are the global solutions to the constrained optimization problem $\min_{\lambda \in K} \varphi^c_{Au}(\lambda,\eta)$ with the parameter $\eta$.
    As a result, for any feasible point that satisfies the constraint (\ref{equation: generalized primal gap constraint based reformulation for OCPEC VI}), $\varphi^c_{Au}(\lambda, \eta) \leq 0$ must be active, and the gradient of $\varphi^c_{Au}$ is either zero or linearly dependent with the gradient of activated $g(\lambda) \geq 0$, which violates LICQ. 
    Similarly, the zeros of $\varphi^{ab}_{Au}$ are the global solutions to the unconstrained optimization problem $\min_{\lambda \in \mathbb{R}^{n_{\lambda}}} \ \varphi^{ab}_{Au}(\lambda,\eta)$, thus $\varphi^{ab}_{Au}(\lambda, \eta) \leq 0$ must be active and its gradient should be zero.
    
    Regarding the MFCQ, it implies the existence of a feasible interior point.
    As has been mentioned, $\varphi^c_{Au}(\lambda, \eta) \leq 0$ must be active for any feasible point satisfying constraints (\ref{equation: generalized primal gap constraint based reformulation for OCPEC VI}).
    Since $\varphi^c_{Au}$ is nonnegative for any $\lambda \in K$, it is impossible to find a point $\lambda \in K$ such that $\varphi^c_{Au}(\lambda, \eta) < 0$ holds, in other words, constraint system (\ref{equation: generalized primal gap constraint based reformulation for OCPEC VI}) does not have a feasible interior and thereby violates MFCQ.
    Similarly, it is impossible to find a point $\lambda \in \mathbb{R}^{n_{\lambda}}$ such that $\varphi^{ab}_{Au}(\lambda, \eta) < 0$ holds, thus constraint system (\ref{equation: generalized D gap constraint based reformulation for OCPEC VI}) also violates MFCQ.
\end{IEEEproof}

\subsection{Proof of Theorems \ref{theorem: exponential convergence of the Y system} and \ref{theorem: error analysis of the Y system (explicit Euler)} }\label{subsection: proof of convergence and error analysis (explicit Euler) of the Y system}

The proof needs properties of the generalized Jacobian.
\begin{proposition}[Proposition 7.1.4, \cite{facchinei2003finite}]\label{proposition: properties of generalized Jacobian}
    Let $G: \Omega \rightarrow \mathbb{R}^m$ be a locally Lipschitz continuous function in an open set $\Omega \subseteq \mathbb{R}^n$. 
    \begin{itemize}
        \item $\partial G (x)$ is \textit{nonempty, convex, and compact} for any $x \in \Omega$;
        \item $\partial G (x)$ is \textit{closed} at $x$, i.e,, for each $\varepsilon > 0$, there is a $\delta > 0$ such that $\partial G(y) \subseteq \partial G(x) + \mathbb{B}(0, \varepsilon), \forall y \in \mathbb{B}(x, \delta)$. 
        \hfill \IEEEQEDopen
    \end{itemize}
\end{proposition}

\begin{lemma}\label{lemma: mean value theorem for Y and Y star}
    Let Assumption \ref{assumption: functions of general parameterized NLP} holds.
    Let $\boldsymbol{Y}(\tau)$ and $s(\tau)$ be the solutions to (\ref{equation: stable system to govern solution trajectory}) and (\ref{equation: stable system to govern parameter trajectory}), respectively.
    For each $\tau \geq 0$, there exists $n_Y$ points $z^i_{\tau}$ in $(\boldsymbol{Y}(\tau), \boldsymbol{Y}^*(\tau))$ and $n_Y$ scalars $\alpha^i_{\tau} \geq 0$ with $\sum_{i = 1}^{n_Y} \alpha^i_{\tau} = 1$ such that
    \begin{equation*}
        \boldsymbol{T}(\boldsymbol{Y}(\tau), s(\tau)) = \boldsymbol{T}(\boldsymbol{Y}^*(\tau), s(\tau)) + \mathcal{M}_{\tau} (\boldsymbol{Y}(\tau) - \boldsymbol{Y}^*(\tau))
    \end{equation*}   
    with $\mathcal{M}_{\tau} = \sum_{i=1}^{n_Y} \alpha^i_{\tau} \mathcal{K}^i_{\tau}$ and $\mathcal{K}^i_{\tau} \in \partial \boldsymbol{T}(z^i_{\tau}, s(\tau))$.
    \hfill \IEEEQEDopen
\end{lemma}
\begin{IEEEproof}
    Since $\boldsymbol{T}(\boldsymbol{Y},s)$ is Lipschitz continuous, this lemma is the direct result of the \textit{mean value theorem} for Lipschitz continuous functions (Proposition 7.1.16, \cite{facchinei2003finite}).
\end{IEEEproof}

We formally state the proof of Theorem \ref{theorem: exponential convergence of the Y system} as follows.
\begin{IEEEproof}  
    We first prove the asymptotic convergence using the candidate Lyapunov function $V(\boldsymbol{Y}, s) = \frac{1}{2} \| \boldsymbol{T}(\boldsymbol{Y}, s) \|_2^2$.
    We have that $V(\boldsymbol{Y}, s) \geq 0$, and $V(\boldsymbol{Y}, s) = 0$ if and only if $\boldsymbol{Y}(\tau) = \boldsymbol{Y}^*(\tau)$.
    The time derivative of $V$ can be written as:
    \begin{equation*}\label{equation: time derivative of Lyapunov function for Y system}
            \Dot{V} = \boldsymbol{T}^T(\mathcal{K} (-\mathcal{K}^{-1}(\epsilon_T \boldsymbol{T} + \mathcal{S} \dot{s} )) + \mathcal{S} \dot{s})  = - 2 \epsilon_T V.       
    \end{equation*}
    Thus, $\Dot{V} < 0$ for all $\boldsymbol{Y}(\tau) \neq \boldsymbol{Y}^*(\tau)$.
    Consequently, following from Theorem 3.3 in \cite{Hessan2015nonlinearcontrol}, there exists a neighborhood of $\boldsymbol{Y}^*_0$ denoted by 
    $\mathcal{N}^{*}_{asy}$, such that for any $\boldsymbol{Y}_0 \in \mathcal{N}^{*}_{asy}$, we have that $\boldsymbol{Y}(\tau)$ asymptotically converges to $\boldsymbol{Y}^*(\tau)$ as $\tau \rightarrow \infty$.
    
    In the following, we prove the exponential convergence, which is inspired by Proposition 2 in \cite{fazlyab2017prediction}.
    First, since $\boldsymbol{Y}(\tau)$ is derived from the stable system (\ref{equation: differential equation for the KKT system (stabilized KKT residual form)}),
    the following inequality holds with a constant $\alpha_T$ satisfying $0 < \alpha_T < \epsilon_T$:
    \begin{equation}\label{equation: bounds for solution of KKT stable system}
        \| \boldsymbol{T}(\boldsymbol{Y}(\tau), s(\tau)) \|_2 \leq \| \boldsymbol{T}(\boldsymbol{Y}(0), s(0)) \|_2 e^{- \alpha_T \tau}.
    \end{equation} 
    
    \textcolor{blue}{Next, we establish the nonsingularity of $\mathcal{M}_{\tau}$ in Lemma \ref{lemma: mean value theorem for Y and Y star}.
    Note that even though we prove in Lemma \ref{lemma: nonsingularity of the KKT matrix} that each $\mathcal{K}^i_{\tau}$ is nonsingular, this does not guarantee that their convex combination $\mathcal{M}_{\tau} = \sum_{i=1}^{n_Y} \alpha^i_{\tau} \mathcal{K}^i_{\tau}$ is also nonsingular.
    To establish the nonsingularity of $\mathcal{M}_{\tau}$, we exploit several properties of the generalized Jacobian, namely closeness and convexity, as presented in Proposition \ref{proposition: properties of generalized Jacobian}}.
    Specifically, based on the closeness of $\partial \boldsymbol{T}$, for each $\tau \geq 0$, we can find a neighborhood of $\partial \boldsymbol{T}(\boldsymbol{Y}^*(\tau), s(\tau))$ defined by $\mathcal{N}^{\varepsilon}_{\tau} := \partial \boldsymbol{T}(\boldsymbol{Y}^*(\tau), s(\tau)) + \mathbb{B}(0, \varepsilon_{\tau})$ with $\varepsilon_{\tau} > 0$, such that $ \partial \boldsymbol{T}(z^i_{\tau}, s(\tau)) \subseteq \mathcal{N}^{\varepsilon}_{\tau}$ for all $z^i_{\tau}$ in $(\boldsymbol{Y}(\tau), \boldsymbol{Y}^*(\tau))$.
    Therefore, $\mathcal{M}_{\tau} $ also belongs to $ \mathcal{N}^{\varepsilon}_{\tau}$ because it is a convex combination of $\mathcal{K}^i_{\tau} \in \partial \boldsymbol{T}(z^i_{\tau}, s(\tau))$.
    Moreover, since $\boldsymbol{Y}(\tau)$ asymptotically converges to $\boldsymbol{Y}^*(\tau)$ as $\tau \rightarrow \infty$, we have that $\{z^i_{\tau}\}_{\tau = 0}^{\infty} \rightarrow \boldsymbol{Y}^*(\tau)$ as $\tau \rightarrow \infty$ for each $i \in \{1, \cdots, n_Y \}$.
    Thus, $\{\varepsilon_{\tau}\}_{\tau = 0}^{\infty} \rightarrow 0$ as $\tau \rightarrow \infty$ and $\{\mathcal{M}_{\tau}\}_{\tau = 0}^{\infty}$ converges to one element in $ \partial \boldsymbol{T}(\boldsymbol{Y}^*(\tau), s(\tau))$, which implies that $\mathcal{M}_{\tau}$ becomes nonsingular as $\tau \rightarrow \infty$.    
    
    Finally, following from Lemma \ref{lemma: mean value theorem for Y and Y star} and (\ref{equation: bounds for solution of KKT stable system}), and the nonsingularity of $\mathcal{M}_{\tau} $, we have:
    \begin{equation*}
        \begin{split}
            &\| \boldsymbol{Y}(\tau) - \boldsymbol{Y}^*(\tau) \|_2 \\
          = & \|\mathcal{M}^{-1}_{\tau} (\boldsymbol{T}(\boldsymbol{Y}(\tau), s(\tau)) - \boldsymbol{T}(\boldsymbol{Y}^*(\tau), s(\tau)))  \|_2 \\
          \leq & \beta_M \| \boldsymbol{T}(\boldsymbol{Y}(0), s(0)) \|_2 e^{- \alpha_T \tau} \\
          \leq & \beta_M L_T \|\boldsymbol{Y}(0) - \boldsymbol{Y}^*(0) \|_2 e^{- \alpha_T \tau},
        \end{split}        
    \end{equation*}
    where $L_T > 0$ is the Lipschitz constant for $\boldsymbol{T}$, and $\beta_M > 0$ is the constant that $ \beta_M \geq \| \mathcal{M}^{-1}_{\tau} \|_2$.
    Thus, the proof is completed with $k_1 = \beta_M L_T$ and $k_2 = \alpha_T$.
\end{IEEEproof}

We formally state the proof of Theorem \ref{theorem: error analysis of the Y system (explicit Euler)} as follows.
\begin{IEEEproof}
    From $\boldsymbol{Y}_{l+1} = \boldsymbol{Y}_l + \Delta \tau \dot{\boldsymbol{Y}}_l$, we first have:
    \begin{equation}\label{equation: error term of explicit Euler method for solution system}
        \boldsymbol{Y}_{l+1} - \boldsymbol{Y}_{l+1}^* = (\boldsymbol{Y}_l - \boldsymbol{Y}_l^*) + \Delta \tau \dot{\boldsymbol{Y}}_l + (\boldsymbol{Y}_l^* - \boldsymbol{Y}_{l+1}^*).
    \end{equation}
    Regarding the term $\Delta \tau\dot{\boldsymbol{Y}}_l$ in (\ref{equation: error term of explicit Euler method for solution system}), it can be written as:
    \begin{equation}\label{equation: expansion of dot_Y_l}
        \begin{split}
            \Delta \tau\dot{\boldsymbol{Y}}_l & = -\Delta \tau\mathcal{K}^{-1}_l(\epsilon_T \boldsymbol{T}_l - \epsilon_s \mathcal{S} (s_l - s_e)) \\
            & = - \epsilon_T \Delta \tau (\boldsymbol{Y}_l - \boldsymbol{Y}^*_l) + \epsilon_s \Delta \tau \mathcal{K}^{-1}_l \mathcal{S} (s_l - s_e),
        \end{split}
    \end{equation}
    where $\boldsymbol{T}_l$ and $\mathcal{K}_l$ are the values of (\ref{equation: KKT condition (equation system) of the general parameterized NLP}) and (\ref{equation: KKT matrix}) with $Y_l$ and $s_l$, respectively.
    The last equality in (\ref{equation: expansion of dot_Y_l}) uses the semismooth Newton approximation of $\boldsymbol{T}(\boldsymbol{Y}, s_l)$ at $\boldsymbol{Y}_l$, i.e., $\boldsymbol{T}(\boldsymbol{Y}^*_l, s_l)  - \boldsymbol{T}(\boldsymbol{Y}_l, s_l) = \mathcal{K}(\boldsymbol{Y}_l, s_l) (\boldsymbol{Y}^*_l - \boldsymbol{Y}_l)$ with $\boldsymbol{T}(\boldsymbol{Y}^*_l, s_l) = 0$.
    Regarding the term $(\boldsymbol{Y}_l^* - \boldsymbol{Y}_{l+1}^*)$ in (\ref{equation: error term of explicit Euler method for solution system}), we have
    \begin{equation}\label{equation: inequality for the change of optimal solution}
            \| \boldsymbol{Y}_{l+1}^* - \boldsymbol{Y}_{l}^*\|_2  \leq L_Y |s_{l+1} - s_{l} | \leq L_Y \beta_s | \dot{s}_l |
    \end{equation}
    with constants $L_Y, \beta_s > 0$ and $\dot{s}_l$ is the value of (\ref{equation: stable system to govern parameter trajectory}) at $\tau_l$.
    The first inequality in (\ref{equation: inequality for the change of optimal solution}) follows from the \textit{implicit function theorem} (Proposition 7.1.18, \cite{facchinei2003finite}) for the Lipschitz continuous equation $\boldsymbol{T}(\boldsymbol{Y}, s) = 0$.
    The second inequality in (\ref{equation: inequality for the change of optimal solution}) holds because (\ref{equation: stable system to govern parameter trajectory}) implies that $|s_{l+1} - s_{l} |$ is always over-predicted by $\Delta \tau | \dot{s}_l |$ (i.e., $|s_{l+1} - s_{l} | < \Delta \tau | \dot{s}_l |$).  
    By substituting (\ref{equation: expansion of dot_Y_l}) into (\ref{equation: error term of explicit Euler method for solution system}), taking the norm inequality, and using (\ref{equation: inequality for the change of optimal solution}), we have:
    \begin{equation*}
        \begin{split}
            & \| \boldsymbol{Y}_{l+1} - \boldsymbol{Y}_{l+1}^* \|_2 \\ 
            \leq & | 1 - \epsilon_T \Delta \tau| \| \boldsymbol{Y}_l - \boldsymbol{Y}^*_l \|_2 + \epsilon_s (\Delta \tau \beta_{KS} + L_Y\beta_s)(s_l - s_e),
        \end{split}
    \end{equation*}
    where $\beta_{KS} > 0$ is a constant that $\beta_{KS} \geq \| \mathcal{K}^{-1}_l \mathcal{S}\|_2$.
    Thus, the proof is completed with $k_3 = \Delta \tau \beta_{KS} + L_Y\beta_s$.
\end{IEEEproof}

\subsection{Proof of Lemma \ref{lemma: nonsingularity of the KKT matrix}}\label{subsection: proof of nonsingularity of the KKT matrix}

\begin{IEEEproof}
    We prove this lemma by contradiction.
    Suppose that $\mathcal{K}$ is singular, then there exists a non-zero vector $q \in \mathbb{R}^{n_Y}$ such that $\mathcal{K} q = 0$.
    By dividing $q = [q^T_1, q^T_2, q^T_3, q^T_4]^T$ with $q_1 \in \mathbb{R}^{n_z}$, $q_2 \in \mathbb{R}^{n_c}$, $q_3 \in \mathbb{R}^{n_h}$, and $q_4 \in \mathbb{R}^{n_c}$, we obtain
    \begin{subequations}\label{equation: KKT matrix times vector}
        \begin{align}
            & \mathcal{H} q_1 + \nabla_{\boldsymbol{z}} \boldsymbol{h}^T  q_3 - \nabla_{\boldsymbol{z}} \boldsymbol{c}^T q_4 = 0, \label{equation: KKT matrix times vector 1} \\
            & \nabla_{\boldsymbol{z}} \boldsymbol{h} q_1 = 0, \label{equation: KKT matrix times vector 2} \\
            & \nabla_{\boldsymbol{z}} \boldsymbol{c} q_1 - q_2 = 0, \label{equation: KKT matrix times vector 3} \\
            & \nabla_{\boldsymbol{v}_c} \Psi q_2 + \nabla_{\boldsymbol{\gamma}_c} \Psi q_4 = 0, \label{equation: KKT matrix times vector 4}
        \end{align}
    \end{subequations}
    Since the strict complementarity condition holds, we have $\nabla_{\boldsymbol{v}_c} \Psi \prec 0 $ and $ \nabla_{\boldsymbol{\gamma}_c} \Psi \prec 0$.
    By substituting (\ref{equation: KKT matrix times vector 3}) and (\ref{equation: KKT matrix times vector 4}) into (\ref{equation: KKT matrix times vector 1}) to eliminate $q_2$ and $q_4$, (\ref{equation: KKT matrix times vector}) becomes
    \begin{equation}\label{equation: simplified KKT matrix times vector}
        \begin{bmatrix}
            \mathcal{H} + \mathcal{R}_c & \nabla_{\boldsymbol{z}}\boldsymbol{h}^T  \\
               \nabla_{\boldsymbol{z}}\boldsymbol{h} & 0
        \end{bmatrix}
        \begin{bmatrix}
            q_1 \\
            q_3
        \end{bmatrix}
        = 0,        
    \end{equation}
    with matrix $\mathcal{R}_c = \nabla_{\boldsymbol{z}} \boldsymbol{c}^T (\nabla_{\boldsymbol{\gamma}_c} \Psi)^{-1} \nabla_{\boldsymbol{v}_c} \Psi  \nabla_{\boldsymbol{z}} \boldsymbol{c} \succ 0$.
    Thus following from Assumption \ref{assumption: CQ recover of general parameterized NLP} and \ref{assumption: reduced Hessian of general parameterized NLP}, the linear system (\ref{equation: simplified KKT matrix times vector}) only has zero solution $q_1 = q_3 = 0$, and thereby $q_2 = q_4 = 0$, which contradicts the assumption made at the beginning that $q$ is a non-zero vector.
    Thus, $\mathcal{K}$ is non-singular.
    Based on Lemma 7.5.2 in \cite{facchinei2003finite}, the nonsingularity holds in the neighborhood of $\boldsymbol{Y}^*$.
\end{IEEEproof}

\subsection{Proof of second and third statements in Proposition \ref{proposition: properties of gap function}}\label{subsection: proof of second and third statements in properties of gap function}

To recap, given a closed convex set $K \subseteq \mathbb{R}^{n_{\lambda}}$ and a vector $\eta \in \mathbb{R}^{n_{\lambda}}$, SOL$(K, \eta)$ is a set consisting of vectors $\lambda \in K$ that satisfy infinitely many inequalities $(\omega - \lambda)^T \eta \geq 0, \forall \omega \in K$.
The proof of the second and third statements in Proposition \ref{proposition: properties of gap function} needs the properties of the saddle problem.
\begin{definition}[Saddle problem]
    Let $X \subseteq \mathbb{R}^n$ and $Y \subseteq \mathbb{R}^m$ be two given closed sets, let $L: X \times Y  \rightarrow \mathbb{R}$ denote an arbitrary function, called a \textit{saddle function}. 
    The saddle problem associated with this triple $(L, X, Y)$ is to find a pair of vectors $(x^*, y^*) \in X \times Y$, called a \textit{saddle point}, such that $L(x^*, y) \leq L(x^*, y^*) \leq L(x,y^*), \forall (x, y) \in X \times Y$. \hfill \IEEEQEDopen
\end{definition}

\begin{proposition}[Theorem 1.4.1, \cite{facchinei2003finite}]\label{proposition: infsup and supinf for saddle problem}
    Let $L : X \times Y \subseteq \mathbb{R}^{n} \times \mathbb{R}^{m} \rightarrow \mathbb{R}$ be a given saddle function. 
    It holds that:
    \begin{equation}\label{equation: relationship bewteen infsup and supinf}
        \inf_{x \in X} \sup_{y \in Y} L(x, y) \geq \sup_{y \in Y} \inf_{x \in X} L(x, y).
    \end{equation}
    Let $\varphi(x) := \sup_{y \in Y} L(x, y)$ and $\psi(y) = \inf_{x \in X} L(x, y)$ be a pair of scalar functions associated with the saddle function $L(x, y)$. 
    Then, for a given pair $(x^*, y^*) \in X \times Y$, the following three statements are equivalent:
    \begin{itemize}
        \item $(x^*, y^*)$ is a saddle point of $L$ on $X \times Y$;
        \item $x^*$ is a minimizer of $\varphi(x)$ on $X$, $y^*$ is a maximizer of $\psi(y)$ on $Y$, and equality holds in (\ref{equation: relationship bewteen infsup and supinf});
        \item $\varphi(x^*) = \psi(y^*) = L(x^*, y^*)$. \hfill \IEEEQEDopen
    \end{itemize}     
\end{proposition}

We now present the proof of the second statement.

\begin{IEEEproof}
    We first prove the nonnegativity property.
    For any given $\hat{\lambda} \in K$ and $\hat{\eta} \in \mathbb{R}^{n_{\lambda}}$, supposing that the maximum of $L^c_{Au}(\hat{\lambda}, \hat{\eta},\omega)$ is obtained at $\hat{\omega} \in K$, then we have $L^c_{Au}(\hat{\lambda}, \hat{\eta}, \hat{\omega}) \geq L^c_{Au}(\hat{\lambda}, \hat{\eta},\omega), \forall \omega \in K$, which includes the case that $\omega = \hat{\lambda}$:
    \begin{equation*}
        \begin{split}
            L^c_{Au}(\hat{\lambda}, \hat{\eta}, \hat{\omega}) & \geq \underbrace{cd(\hat{\lambda}) - cd(\hat{\lambda}) + (\hat{\eta}^T - c\nabla_{\lambda} d(\hat{\lambda}))(\hat{\lambda}-\hat{\lambda})}_{L^c_{Au}(\hat{\lambda}, \hat{\eta}, \hat{\lambda})} \\
            & = 0.
        \end{split}         
    \end{equation*}
    Thus, we have $\varphi^c_{Au}(\lambda, \eta):= \sup_{\omega \in K} L^c_{Au}(\lambda, \eta, \omega) \geq 0, \forall \lambda \in K$, and the nonnegative property is proved.
    Similarly, we also have $\psi^c_{Au}(\eta, \omega) :=\inf_{\lambda \in K} L^c_{Au}(\lambda, \eta, \omega) \leq 0, \forall \omega \in K$.
    
    We next prove the sufficient condition of the equivalence property, i.e., $\varphi^{c}_{Au}(\lambda, \eta) = 0$ with $\lambda \in K \Rightarrow \lambda \in$ SOL$(K, \eta)$. 
    From Proposition \ref{proposition: infsup and supinf for saddle problem} and the properties that $\varphi^c_{Au}(\lambda, \eta) \geq 0$ and $\psi^c_{Au}(\eta,\omega) \leq 0$, for any given $\eta^*$, we have that $\varphi^c_{Au}(\lambda, \eta^*) =0$ if and only if $\lambda = \lambda^*$, where $\lambda^*$ is the primal part of the saddle point $(\lambda^*, \omega^*)$ of $L^c_{Au}(\lambda, \eta^*, \omega)$, that is:
    \begin{equation*}
        \varphi^c_{Au}(\lambda^*, \eta^*) = L^c_{Au}(\lambda^*, \eta^*, \omega^*) = \psi^c_{Au}(\eta^*, \omega^*) = 0.
    \end{equation*}
    From the definition $\varphi^c_{Au}(\lambda^*, \eta^*) = \sup_{\omega \in K} L^c_{Au}(\lambda^*, \eta^*, \omega)$, we have that
    \begin{equation*}
        \begin{split}
            & \underbrace{cd(\lambda^*) - cd(\omega) + ((\eta^*)^T - c\nabla_{\lambda} d(\lambda^*))(\lambda^*-\omega)}_{L^c_{Au}(\lambda^*, \eta^*, \omega)}\\
            \leq \ & \varphi^c_{Au}(\lambda^*, \eta^*) = 0, \quad \forall \omega \in K, \\
        \end{split}         
    \end{equation*}
    and the maximum of $L^c_{Au}(\lambda^*, \eta^*, \omega)$ can be obtained at $\omega = \lambda^*$. 
    Thus, we have the first-order primal necessary condition:
    \begin{equation*}
        \begin{split}
              & \underbrace{ (-c\nabla_{\lambda} d(\lambda^*) - ((\eta^*)^T - c\nabla_{\lambda} d(\lambda^*)) )}_{\nabla_{\omega} L^c_{Au}(\lambda^*, \eta^*, \lambda^*)} (\omega - \lambda^*) \\
            = & -(\eta^*)^T (\omega - \lambda^*) \leq  0, \quad \forall \omega \in K,            
        \end{split}         
    \end{equation*}
    which means that $\lambda^*$ solves the VI$(K, \eta^*)$.    

    We finally prove the necessary condition of the equivalence property, i.e., $\lambda \in$ SOL$(K, \eta) \Rightarrow \varphi^c_{Au}(\lambda, \eta) = 0$.
    For any given $\hat{\eta} \in \mathbb{R}^{n_{\lambda}}$, suppose that $\hat{\lambda} \in$ SOL$(K, \hat{\eta})$, from the definition of SOL$(K, \hat{\eta})$, we have:
    \begin{equation*}
        \underbrace{ (-c\nabla_{\lambda} d(\hat{\lambda}) - \hat{\eta}^T + c\nabla_{\lambda} d(\hat{\lambda}) )}_{\nabla_{\omega} L^c_{Au}(\hat{\lambda}, \hat{\eta}, \hat{\lambda})} (\omega - \hat{\lambda}) \leq  0, \quad \forall \omega \in K.      
    \end{equation*}
    This implies that the maximum of $L^c_{Au}(\hat{\lambda}, \hat{\eta}, \omega)$ can be obtained at $\omega = \hat{\lambda}$, which is $L^c_{Au}(\hat{\lambda}, \hat{\eta}, \hat{\lambda}) = 0$.
    Hence we have $\varphi^c_{Au}(\hat{\lambda}, \hat{\eta}) = 0$.
\end{IEEEproof}

\begin{remark}
    Our proof of the nonnegativity and equivalence properties is based on Proposition \ref{proposition: infsup and supinf for saddle problem} and the optimality conditions in the form of VI, which is slightly different from \cite{auchmuty1989variational}, where Auchmuty proves these properties using Proposition \ref{proposition: infsup and supinf for saddle problem} and generalized Young's inequality. \hfill \IEEEQEDopen
\end{remark}

The proof of the third statement needs the following lemma about the properties of the generalized D-gap function:
\begin{lemma}\label{lemma: properties of generalized D gap function}
    The generalized D-gap function $\varphi^{ab}_{Au}(\lambda, \eta)$ satisfies the following inequalities:
    \begin{equation}\label{equation: inequalities properties of generalized D-gap function}
        \varphi^{ab}_{Au}(\lambda, \eta) \geq \frac{m(b-a)}{2} \|\hat{\omega}^b - \lambda\|^2_2,
    \end{equation}
    with $m > 0$ a constant for the strong convexity of $d$, that is, $d(\omega) \geq d(\lambda) + \nabla_{\lambda} d(\lambda) (\omega - \lambda) + \frac{m}{2} \|\omega - \lambda\|^2_2$.
    \hfill \IEEEQEDopen
\end{lemma}

\begin{IEEEproof}
    The proof is inspired by Lemma 10.3.2 in \cite{facchinei2003finite}.
    The inequality (\ref{equation: inequalities properties of generalized D-gap function}) is derived by:
    \begin{equation*}
        \begin{split}
            & \varphi^{ab}_{Au}(\lambda, \eta) \\
        = &\ \varphi^{a}_{Au}(\lambda, \eta) - \varphi^{b}_{Au}(\lambda, \eta) \\
        \geq &\ \eta^T(\lambda - \hat{\omega}^b) + a(d(\lambda) - d(\hat{\omega}^b) + \nabla_{\lambda} d(\lambda) (\hat{\omega}^b-\lambda)) \\
            & - \eta^T(\lambda - \hat{\omega}^b) - b(d(\lambda) - d(\hat{\omega}^b) + \nabla_{\lambda} d(\lambda) (\hat{\omega}^b-\lambda)) \\
        = &\ -(b-a)(d(\lambda) - d(\hat{\omega}^b) + \nabla_{\lambda} d(\lambda) (\hat{\omega}^b-\lambda))  \\
        \geq & \ \frac{m(b-a)}{2} \|\hat{\omega}^b - \lambda\|^2_2.
        \end{split}         
    \end{equation*}
\end{IEEEproof}

We now present the proof of the third statement.

\begin{IEEEproof}
    The proof is inspired by Theorem 10.3.3 in \cite{facchinei2003finite}.
    Regarding the nonnegativity property that $\varphi^{ab}_{Au}(\lambda, \eta) \geq 0, \forall \lambda \in \mathbb{R}^n$, it follows from (\ref{equation: inequalities properties of generalized D-gap function}).
    For the sufficient condition of the equivalence property, i.e., $\varphi^{ab}_{Au}(\lambda, \eta) = 0 \Rightarrow \lambda \in$ SOL$(K, \eta)$, if $\varphi^{ab}_{Au}(\lambda, \eta) = 0$, from (\ref{equation: inequalities properties of generalized D-gap function}) we have $\lambda = \hat{\omega}^b$,
    which implies that $\lambda \in K$ and $\varphi^{b}_{Au}(\lambda, \eta) = 0$, hence $\lambda \in$ SOL$(K, \eta)$.
    For the necessary condition of the equivalence property, i.e., $\lambda \in$ SOL$(K, \eta) \Rightarrow \varphi^{ab}_{Au}(\lambda, \eta) = 0$, since $\lambda \in$ SOL$(K, \eta)$, based on the equivalence property of $\varphi^{c}_{Au}(\lambda, \eta)$, we have $\varphi^{a}_{Au}(\lambda, \eta) = \varphi^{b}_{Au}(\lambda, \eta) = 0$  hence $\varphi^{ab}_{Au}(\lambda, \eta) = 0$.
\end{IEEEproof}

\bibliographystyle{IEEEtran}
\bibliography{IEEEabrv, reference}

\end{document}